\documentclass[a4paper,11pt]{article}
\usepackage{aligned-overset}
\usepackage{amsmath,amsthm}
\usepackage{authblk}
\usepackage[style=numeric-comp,maxbibnames=99,bibencoding=utf8,firstinits=true]{biblatex}
\usepackage{booktabs}
\usepackage{cancel}
\usepackage{cases}
\usepackage{colortbl}
\usepackage[hmargin={26mm,26mm},vmargin={30mm,35mm}]{geometry}
\usepackage{nicefrac}
\usepackage[colorlinks,allcolors={blue}]{hyperref}
\usepackage{paralist}
\usepackage{enumitem}
\usepackage[compact,small]{titlesec}
\usepackage{tikz-cd}
\usepackage{graphicx}
\usepackage{microtype}
\usepackage{placeins}
\usepackage{subcaption}
\usepackage{newtxtext}
\usepackage{newtxmath}
\usepackage{xspace}

\bibliography{ddr-sobolev-ns}

\newcommand{\email}[1]{\href{mailto:#1}{#1}}

\numberwithin{equation}{section}

\newtheorem{theorem}{Theorem}
\newtheorem{proposition}[theorem]{Proposition}
\newtheorem{lemma}[theorem]{Lemma}

\theoremstyle{remark}
\newtheorem{remark}[theorem]{Remark}
\theoremstyle{definition}

\newcommand{\Real}{\mathbb{R}}

\DeclareRobustCommand{\bvec}[1]{\boldsymbol{#1}}
\newcommand{\uvec}[1]{\underline{\bvec{#1}}}
\newcommand{\cvec}[1]{\bvec{\mathcal{#1}}}

\DeclareMathOperator{\GRAD}{\bf grad}
\DeclareMathOperator{\CURL}{\bf curl}
\DeclareMathOperator{\DIV}{div}

\DeclareMathOperator{\VROT}{\bf rot}

\newcommand{\compl}{{\rm c}}

\newcommand{\Hgrad}[1]{H(\GRAD;#1)}
\newcommand{\Hgradz}[1]{H_0(\GRAD;#1)}
\newcommand{\Hcurl}[1]{\bvec{H}(\CURL;#1)}
\newcommand{\Hcurlz}[1]{\bvec{H}_0(\CURL;#1)}

\newcommand{\Hdiv}[1]{\bvec{H}(\DIV;#1)}
\newcommand{\Hdivz}[1]{\bvec{H}_0(\DIV;#1)}
\newcommand{\Hdivzz}[1]{\bvec{H}_0(\DIV=0;#1)}
\newcommand{\Hbullet}[1]{H(\bullet;#1)}

\newcommand{\Xgrad}[2]{{X}_{\GRAD}^{#1}(#2)}
\newcommand{\Xcurl}[2]{{\bvec{X}}_{\CURL}^{#1}(#2)}
\newcommand{\Xdiv}[2]{{\bvec{X}}_{\DIV}^{#1}(#2)}
\newcommand{\Xbullet}[2]{{X}_{\,\bullet}^{#1}(#2)}

\newcommand{\Igrad}[2]{\underline{I}_{\GRAD,#2}^{#1}}
\newcommand{\qIgrad}[2]{\widehat{\underline{I}}_{\GRAD,#2}^{#1}}
\newcommand{\Icurl}[2]{\uvec{I}_{\CURL,#2}^{#1}}
\newcommand{\qIcurl}[2]{\widehat{\uvec{I}}_{\CURL,#2}^{#1}}
\newcommand{\wIcurl}[2]{\widetilde{\uvec{I}}_{#2}}
\newcommand{\qIdiv}[2]{\widehat{\uvec{I}}_{\DIV,#2}^{#1}}
\newcommand{\qIbullet}[2]{\widehat{\underline{I}}_{\bullet,#2}^{#1}}
\newcommand{\Idiv}[2]{\uvec{I}_{\DIV,#2}^{#1}}

\newcommand{\Liftgrad}{{\mathcal L}_{\GRAD,h}}
\newcommand{\Liftcurl}{\bvec{\mathcal L}_{\CURL,h}}
\newcommand{\Liftdiv}{\bvec{\mathcal L}_{\DIV,h}}
\newcommand{\Liftbullet}{{\mathcal L}_{\bullet,h}}

\newcommand{\NEsymbol}{\boldsymbol{\mathcal{N}}}
\newcommand{\RTsymbol}{\boldsymbol{\mathcal{RT}}}
\newcommand{\NE}[1]{\NEsymbol^{#1}}
\newcommand{\RT}[1]{\RTsymbol^{#1}}
\newcommand{\NEproj}[2]{\bvec{\symbolproj}_{\NEsymbol,#2}^{#1}}
\newcommand{\RTproj}[2]{\bvec{\symbolproj}_{\RTsymbol,#2}^{#1}}

\newcommand{\symbolproj}{\pi}

\newcommand{\lproj}[2]{\symbolproj_{\mathcal{P},#2}^{#1}}
\newcommand{\vlproj}[2]{\boldsymbol{\symbolproj}_{\cvec{P},#2}^{#1}}

\newcommand{\uG}[2]{\uvec{G}_{#2}^{#1}}
\newcommand{\uC}[2]{\uvec{C}_{#2}^{#1}}
\newcommand{\D}[2]{D_{#2}^{#1}}

\newcommand{\trFt}{\bvec{\gamma}_{{\rm t},F}^k}

\newcommand{\cells}[1]{\mathcal{T}_{#1}}
\newcommand{\faces}[1]{\mathcal{F}_{#1}}
\newcommand{\edges}[1]{\mathcal{E}_{#1}}
\newcommand{\vertices}[1]{\mathcal{V}_{#1}}

\newcommand{\FT}{\faces{T}}
\newcommand{\ET}{\edges{T}}
\newcommand{\EF}{\edges{F}}
\newcommand{\VT}{\vertices{T}}

\newcommand{\normal}{\bvec{n}}
\newcommand{\tangent}{\bvec{t}}

\newcommand{\Poly}[2][]{\mathcal{P}_{#1}^{#2}}

\newcommand{\vPoly}[2][]{\cvec{P}_{#1}^{#2}}
\newcommand{\Roly}[1]{\cvec{R}^{#1}}
\newcommand{\Goly}[1]{\cvec{G}^{#1}}
\newcommand{\cRoly}[1]{\cvec{R}^{\compl,#1}}
\newcommand{\cGoly}[1]{\cvec{G}^{\compl,#1}}

\newcommand{\vC}[1]{\boldsymbol{C}^{#1}}

\newcommand{\norm}[2]{\|#2\|_{#1}}
\newcommand{\seminorm}[2]{|#2|_{#1}}
\newcommand{\vvvert}{\vert\kern-0.25ex\vert\kern-0.25ex\vert}
\newcommand{\tnorm}[2]{\vvvert #2\vvvert_{#1}}

\newcommand{\term}{\mathfrak{T}}

\DeclareMathOperator{\Image}{Im}

\newcommand{\Mh}{\mathcal{M}_{h}}
\newcommand{\Th}{\mathcal{T}_{h}}
\newcommand{\Fh}{\mathcal{F}_{h}}
\newcommand{\Eh}{\mathcal{E}_{h}}
\newcommand{\Vh}{\mathcal{V}_h}

\newcommand{\Pgrad}[2]{P_{\GRAD,#2}^{#1}}
\newcommand{\Pcurl}[2]{\bvec{P}_{\CURL,#2}^{#1}}
\newcommand{\Pdiv}[2]{\bvec{P}_{\DIV,#2}^{#1}}
\newcommand{\Pbullet}[2]{P_{\bullet,#2}^{#1}}

\newcommand{\Cpoin}{C_{{\rm p},\CURL}}

\newcommand{\Reynolds}{\mathrm{Re}}

\newcommand{\hardcore}{\texttt{HArDCore}\xspace}
\newcommand{\polyrust}{\texttt{poly\_rust}\xspace}

\begin{document}

\title{Sobolev inequalities for the Discrete de Rham complex and an application to the Navier--Stokes equations}

\author[1]{Daniele A. Di Pietro}
\author[1,2]{J\'er\^ome Droniou}
\author[3]{Marien Hanot}
\author[1]{Silvano Pitassi}
\affil[1]{IMAG, Univ Montpellier, CNRS, Montpellier, France, \email{daniele.di-pietro@umontpellier.fr}, \email{jerome.droniou@cnrs.fr}, \email{silvano.pitassi@umontpellier.fr}}
\affil[2]{School of Mathematics, Monash University, Melbourne, Australia}
\affil[3]{University of Lille, Lille, France, \email{marien-lorenzo.hanot@univ-lille.fr}}

\maketitle

\begin{abstract}
  We prove Sobolev inequalities for the discrete de Rham (DDR) method, an arbitrary-order discretisation of the de Rham complex on polyhedral meshes.
  We show that, for any polyhedral domain $\Omega$, the DDR counterparts of the spaces $\Hgrad{\Omega}$, $\Hcurl{\Omega}$ and $\Hdiv{\Omega}$ satisfy the same Sobolev embeddings as their continuous counterparts. For the curl and divergence spaces, this establishes on quasi-uniform meshes an $L^r$-bound (with $r>2$ in general, and $r=6$ if $\Omega$ is convex) on discrete vector fields that are orthogonal to discrete gradients (resp.~curls), in terms of the $L^2$-norms of the vector fields themselves and their discrete curl (resp.~divergence).
  These Sobolev embeddings are the first such results for a polytopal complex. Their proofs hinge on recently introduced quasi-interpolators and conforming liftings.
  As an application, we use these estimates to prove the convergence of a DDR scheme for the Navier--Stokes equations in the curl--curl formulation. The convergence proof uses a Maxwell compactness theorem for DDR and requires no smallness assumptions on the data.
  \medskip\\
  \textbf{Key words.}
  discrete Sobolev inequalities, %
  de Rham complex, %
  discrete de Rham method, %
  polyhedral methods, %
  compatible discretisations, %
  incompressible Navier--Stokes problem. %
  \medskip\\
  \textbf{MSC2020.} 65N12, 
  65N30, 
  14F40, 
  76D05  
\end{abstract}



\section{Introduction}

In this paper, we prove Sobolev inequalities for the polytopal Discrete de Rham (DDR) complex of
\cite{Di-Pietro.Droniou:23}. The proof exploits recent advances in the analysis of this complex, including
conforming liftings \cite{Di-Pietro.Droniou.ea:25} and quasi-interpolators \cite{Chaumont-Frelet.Droniou.ea:25}, to prove $L^r$-bounds on DDR vectors in terms of $L^2$-graph norms that involve discrete counterparts of the gradient, curl or divergence.
These inequalities fit into the general framework of discrete functional analysis and, for the curl and divergence, are the first of their kind for a polytopal complex.
We then exploit these Sobolev inequalities to rigorously prove the convergence of the DDR scheme for the Navier--Stokes problem of \cite{Di-Pietro.Droniou.ea:24} without smallness assumptions on the data, using the discrete Maxwell compactness result of \cite{Chaumont-Frelet.Droniou.ea:25}.
\smallskip

Discrete functional analysis (DFA) is a branch of numerical analysis that establishes discrete counterparts of functional analysis results commonly used in the study of partial differential equations (PDEs).
Examples include Poincaré and Sobolev inequalities, as well as compactness results of Rellich or Aubin--Simon type.
These results are essential for carrying out rigorous convergence analyses of numerical schemes for complex PDEs in real-world settings, e.g., with irregular coefficients and nonlinearities.
Examples of problems in which accounting for these features is crucial include flows in heterogeneous porous media and the Navier--Stokes equations.

For conforming discretisations of coercive linear problems, DFA results follow directly from continuous functional analysis theorems.
For non-conforming methods, on the other hand, establishing these
results is more challenging and often requires advanced technical arguments. Pioneering results for classical finite volume methods can be found in \cite{Eymard.Gallouet.ea:00} and have since been extended to more modern methods; see the review \cite{Droniou:14} for details and a range of complex models for which finite volume methods have been analysed using DFA tools. Discrete inequalities and compactness for discontinuous Galerkin methods are covered, e.g., in \cite{Di-Pietro.Ern:10,Riviere.Walkington:11}.
Efforts have also been made to develop a unified presentation of numerical schemes for classes of linear or nonlinear problems. A particularly relevant example of a unifying paradigm is provided by
the Gradient Discretisation Method, a framework in which generic DFA theorems covering a range of methods have been established \cite{Droniou.Eymard.ea:18}; these include inequalities and compactness results for both stationary and time-dependent problems in various topologies. More recent studies continue to explore discrete inequalities, in particular to broaden their scope of application and/or improve the corresponding constants \cite{Botti.Mascotto:26}.
\smallskip

The aforementioned works, however, only cover DFA for schemes in which the gradient is the main operator. Many important physical models -- in particular those involving electromagnetic phenomena --
also involve the curl and divergence operators, usually organised in the setting of the de Rham complex.
Even at the continuous level, Sobolev inequalities for such operators are much more difficult to obtain than those for the gradient, because curl and divergence have large kernels and do not encode all partial derivatives. Standard approaches deduce these inequalities from $H^s(\Omega)$-regularity for $\Hcurl{\Omega}$ (resp.~$\Hdiv{\Omega}$) vector fields that are $L^2$-orthogonal to all gradients (resp.~curls).
This regularity follows from fine estimates on Hodge decompositions (i.e., ultimately, on solutions to PDEs); see, e.g., \cite{Amrouche.Bernardi.ea:98} and references therein.
Transferring Sobolev inequalities to the discrete setting brings another important challenge: in view of applications to numerical schemes, the orthogonality that can be expected at the discrete level is with respect to gradients (resp.~curls) of discrete functions, not of all functions in $\Hgrad{\Omega}$ (resp.~$\Hcurl{\Omega}$). As a consequence, even for conforming finite elements, very few discrete Sobolev inequalities seem to be available; the main reference we are aware of is \cite{Girault:90}, which establishes these inequalities for Nédélec--Raviart--Thomas finite elements.
\smallskip

Finite element methods are well established and robust, but impose restrictions on the meshes they support.
These limitations prevent, e.g., the use of local adaptation strategies such as non-conforming mesh refinement or coarsening,
and may require more elements than necessary to mesh complex geometries.
The past two decades have therefore seen a push towards polytopal methods, which support meshes made of general polygonal or polyhedral elements; see, e.g., \cite{Beirao-da-Veiga.Brezzi.ea:13,Di-Pietro.Ern.ea:14,Beirao-da-Veiga.Brezzi.ea:16,Di-Pietro.Droniou.ea:20}. Some of these methods have been demonstrated to match or improve on the accuracy of finite element methods at a (sometimes vastly) reduced cost \cite{Dassi.Lovadina.ea:21,Di-Pietro.Droniou:22}. Gradient-based polytopal methods have been studied for some time, and a large body of DFA results has therefore already been established for them (see, for example, \cite{Di-Pietro.Ern:12,Droniou.Eymard.ea:18,Di-Pietro.Droniou:20}).
On the other hand, DFA results for polytopal methods involving discrete counterparts of the curl or divergence are scarcer and much more recent, and mostly concern Poincaré inequalities (sometimes implicit through stability estimates); see \cite{Beirao-da-Veiga.Mascotto.ea:22,Di-Pietro.Droniou:23,Di-Pietro.Hanot:24,Pitassi:26}. Some recent results cover the Maxwell compactness of such methods \cite{Lemaire.Pitassi:25,Chaumont-Frelet.Droniou.ea:25} but, to our knowledge, no discrete Sobolev inequalities have been established for curl/div polytopal methods so far.
\smallskip

The purpose of this article is to fill this gap and contribute to this long line of DFA results for numerical methods. We consider the Discrete de Rham (DDR) method in 3D, an arbitrary-order discretisation of the de Rham complex that is applicable to general polyhedral meshes, and prove Sobolev inequalities for all the spaces involved.
The proof relies on two main ingredients: the \emph{conforming liftings} of \cite{Di-Pietro.Droniou.ea:25}, one-to-one mappings with suitable bounds between the DDR spaces and their continuous counterparts,
and the \emph{commuting quasi-interpolators} of \cite{Chaumont-Frelet.Droniou.ea:25}, which map each conforming function with minimal regularity to a vector in the corresponding DDR space.
The fact that the quasi-interpolator is a left inverse of the lifting
is the key to transferring the continuous Hodge decomposition to the discrete setting, and to exploiting the discrete orthogonality assumption on which the Sobolev embedding relies.
We then combine the newly established discrete Sobolev embedding with the Maxwell compactness result of \cite{Chaumont-Frelet.Droniou.ea:25} to perform the convergence analysis of the DDR scheme for the Navier--Stokes equations, originally proposed in \cite{Di-Pietro.Droniou.ea:24}, without smallness assumptions on the data.

Even though we focus here on the DDR method, our approach is generic and can in all likelihood be extended to virtual element method (VEM) complexes \cite{Beirao-da-Veiga.Brezzi.ea:16}. Indeed, as seen in \cite{Beirao-da-Veiga.Dassi.ea:22}, the lowest-order, lowest-dimensional degrees of freedom (DOFs) in each VEM and DDR space coincide, and since the conforming liftings and quasi-interpolators of \cite{Di-Pietro.Droniou.ea:25,Chaumont-Frelet.Droniou.ea:25} are based precisely on these DOFs and the continuity they represent, designing these tools for VEM should be feasible. Once this is done, the approach we present here can be applied directly to obtain discrete Sobolev inequalities for the VEM complex.
\smallskip

The paper is organised as follows. In the next section, we introduce notation for the continuous function spaces $\Hgrad{\Omega}$, $\Hcurl{\Omega}$ and $\Hdiv{\Omega}$, and briefly recall the related Sobolev inequalities. Section \ref{sec:ddr} presents the Discrete de Rham complex, focusing on its spaces and the relevant norms; full formulae for the discrete operators are not presented here as they are not directly relevant to our purpose. In Section \ref{sec:sobolev}, we state and prove the Sobolev inequalities for all three DDR spaces, discrete versions of $\Hgrad{\Omega}$, $\Hcurl{\Omega}$ and $\Hdiv{\Omega}$. In Section \ref{sec:ddr.ns}, we present the DDR scheme for the Navier--Stokes equations in the curl--curl formulation, complete the error estimate from \cite{Di-Pietro.Droniou.ea:24} (which was proved there under the \emph{assumption} that a discrete Sobolev embedding held), and provide a new convergence proof for that scheme, based on compactness techniques and valid without any smallness assumption on the data. Finally, Section \ref{sec:num.tests} provides some numerical illustrations of the scheme.

\section{Continuous spaces and inequalities}\label{sec:continuous}

In the following, $\Omega$ is a bounded polyhedral domain in $\Real^3$. We use the following notation for the classical Sobolev spaces associated with the gradient, curl and divergence operators:
\begin{align*}
  \Hgrad{\Omega}&\coloneq \{q\in L^2(\Omega)\,:\,\GRAD q\in \bvec{L}^2(\Omega)\},\\
  \Hcurl{\Omega}&\coloneq \{\bvec{v}\in \bvec{L}^2(\Omega)\,:\,\CURL\bvec{v}\in \bvec{L}^2(\Omega)\},\\
  \Hdiv{\Omega}&\coloneq \{\bvec{w}\in \bvec{L}^2(\Omega)\,:\,\DIV\bvec{w}\in L^2(\Omega)\},
\end{align*}
where $\bvec{L}^2(\Omega)\coloneq L^2(\Omega)^3$. We denote by $\Hcurlz{\Omega}$ the subspace of $\Hcurl{\Omega}$ with zero tangential boundary trace, and by $\Hdivz{\Omega}$ the subspace of $\Hdiv{\Omega}$ with zero normal boundary trace.

The classical Sobolev embedding
\begin{equation}\label{eq:continuous.sobolev.grad}
  \Hgrad{\Omega}\hookrightarrow L^6(\Omega)
\end{equation}
holds for any bounded Lipschitz domain \cite[Corollary 9.14]{Brezis:11} (and hence, in particular, for any bounded polyhedral domain).
Sobolev embeddings involving the curl and divergence are less straightforward to state (and prove), as they depend on the geometry of the domain. Specifically, for polyhedral domains, \cite[Proposition 3.7]{Amrouche.Bernardi.ea:98} ensures that there exists $s\in \big\lparen \frac12,1\big\rbrack$, depending on $\Omega$, such that
\begin{equation}\label{eq:regularity.Hcurl.Hdiv}
  \text{$\Hcurl{\Omega}\cap \Hdivz{\Omega}\hookrightarrow \bvec{H}^{s}(\Omega)$ and
    $\Hcurlz{\Omega}\cap \Hdiv{\Omega}\hookrightarrow \bvec{H}^{s}(\Omega)$,}
\end{equation}
where $\bvec{H}^{s}(\Omega)=H^s(\Omega)^3$.
The continuous Sobolev embedding for $H^s(\Omega)$ then yields, with $s^*\coloneq\frac{6}{3-2s}$,
\begin{equation}\label{eq:continuous.sobolev.curl.div}
  \text{$\Hcurl{\Omega}\cap \Hdivz{\Omega}\hookrightarrow \bvec{L}^{s^*}(\Omega)$ and
    $\Hcurlz{\Omega}\cap \Hdiv{\Omega}\hookrightarrow \bvec{L}^{s^*}(\Omega)$.}
\end{equation}
The exponent $s$ can be estimated in terms of the regularity properties of the Laplacian in $\Omega$; see \cite[Remark 3.8]{Amrouche.Bernardi.ea:98}.
If $\Omega$ is convex, we can take $s=1$ and thus $s^*=6$.
Otherwise, for a generic polyhedral domain, we deduce from \cite[Theorem 8.2.1.2 and Remark 8.2.1.4]{Grisvard:85} that $s \leq \min\{\frac{\pi}{\omega}\,:\,\omega\text{ is a dihedral angle between two faces of $\partial\Omega$}\}$.


\section{Discrete de Rham complex}\label{sec:ddr}

\subsection{Mesh and notation}\label{sec:mesh}

The mesh framework we consider is standard in polytopal methods; here, we simply recall some notation following \cite{Di-Pietro.Droniou:23}.
A mesh of $\Omega$ is a set $\Mh=\Th\cup\Fh\cup\Eh\cup\Vh$ of mesh entities, where $\Th$ collects the open disjoint polytopal elements, $\Fh$ is the set of open planar polygonal faces, $\Eh$ is the set of open straight edges and $\Vh$ the set of vertices. We have $\overline{\Omega} = \bigcup_{T\in\Th}\overline{T}$ and $h=\max_{T\in\Th}h_T>0$, where $h_Y$ is the diameter of $Y \in \Mh$. For any $Y\in \Mh$ and $\mathcal X\in\{\mathcal T,\mathcal F,\mathcal E,\mathcal V\}$, we denote by $\mathcal X_Y$ the set of entities in $\mathcal X_h$ that are contained in $\partial Y$ (if $Y$ has a higher dimension than these entities) or whose boundary contains $Y$ (otherwise). We assume that, for $d\in\{1,2,3\}$, any $(d-1)$-dimensional mesh entity is contained in the boundary of some $d$-dimensional mesh entity and that, for any $d$-dimensional $Y\in\Mh$, $\partial Y$ is the union of the lower-dimensional mesh entities contained in $\partial Y$.
The orientation of mesh entities of dimension $1$ and $2$ is fixed as follows:
\begin{itemize}[itemsep=0em,leftmargin=1em,topsep=0.5em]
\item For each $F\in\Fh$, a unit normal vector $\normal_F$ is fixed and, if $T\in\cells{F}$, we denote by $\omega_{TF} \in \{-1, 1\}$ the number such that $\omega_{TF} \normal_F$ points out of $T$.
\item For each $E\in\Eh$, a unit tangent vector $\tangent_E$ is fixed and, if $F\in\faces{E}$, we denote by $\normal_{FE}$ the unit normal vector to $E$ in the plane spanned by $F$ such that $(\tangent_E,\normal_{FE},\normal_F)$ is right-handed; moreover, $\omega_{FE} \in \{-1, 1\}$ is such that $\omega_{FE} \normal_{FE}$ points out of $F$.
\end{itemize}

We also assume that the mesh is such that $(\Th,\Fh)$ satisfies the design and regularity conditions in \cite[Definitions 1.4 and 1.9]{Di-Pietro.Droniou:20}.
We use the shorthand $a\lesssim b$ for $a\le Cb$ with $C$ depending only on $\Omega$, the mesh regularity parameter and, when polynomial functions are involved, the corresponding polynomial degree. The notation $a\simeq b$ means ``$a\lesssim b$ and $b\lesssim a$''.
The mesh regularity assumption ensures that, for each element or face $Y\in\Th\cup\Fh$, we can find~$\bvec{x}_Y\in Y$ such that $Y$ contains a ball/disc centred at $\bvec{x}_Y$ and of diameter $\gtrsim h_Y$.

\subsection{Polynomial spaces}

If $\ell \ge 0$ is an integer and $Y\in\Mh$, the space of polynomial functions of total degree $\le \ell$ on $Y$ is denoted by $\Poly{\ell}(Y)$, with the convention that
$\Poly{-1}(Y)\coloneq\left\{0\right\}$. For $T\in\Th$, we set $\vPoly{\ell}(T)\coloneq \Poly{\ell}(T)^3$ while, for $F\in\Fh$, $\vPoly{\ell}(F)$ denotes the space of vector-valued polynomials $F\to\Real^3$ of degree $\le\ell$ that are everywhere tangent to $F$ (this space is isomorphic to $\Poly{\ell}(F)^2$).
The $L^2(Y)$-orthogonal projector on $\Poly{\ell}(Y)$ is $\lproj{\ell}{Y}$, while the projector on $\vPoly{\ell}(Y)$ is $\vlproj{\ell}{Y}$.
The spaces of broken piecewise polynomials on $\Th$ of degree $\le \ell$ are $\Poly{\ell}(\Th)$ (scalar version) and $\vPoly{\ell}(\Th)$ (vector-valued version); the projectors on these spaces are denoted, respectively, by $\lproj{\ell}{h}$ and $\vlproj{\ell}{h}$.

For $F\in \Fh$, we will need the following polynomial subspaces of $\vPoly{\ell}(F)$:
\[
\Roly{\ell}(F) \coloneq  \VROT_F\Poly{\ell+1}(F), \qquad \cRoly{\ell}(F)  \coloneq (\bvec{x}-\bvec{x}_F)\Poly{\ell-1}(F),
\]
where $\VROT_F$ is the rotor on $F$ (obtained by rotating the tangential gradient).
For  $T\in\Th$, the following subspaces of $\vPoly{\ell}(T)$ will also be useful:
\[
\begin{aligned}
  \Goly{\ell}(T) & \coloneq \GRAD\Poly{\ell+1}(T), & \cGoly{\ell}(T) & \coloneq (\bvec{x}-\bvec{x}_T)\times\vPoly{\ell-1}(T), \\
  \Roly{\ell}(T) & \coloneq \CURL\vPoly{\ell+1}(T), & \cRoly{\ell}(T) & \coloneq (\bvec{x}-\bvec{x}_T) \Poly{\ell-1}(T).
\end{aligned}
\]
For $Y\in\Th\cup\Fh$, the Raviart--Thomas and Nédélec spaces of degree $\ell$ are defined (the latter only for $Y\in\Th$) by
\[
\RT{\ell}(Y)\coloneq \vPoly{\ell-1}(Y)+\cRoly{\ell}(Y)\,,\quad\NE{\ell}(Y)\coloneq \vPoly{\ell-1}(Y)+\cGoly{\ell}(Y).
\]
The $L^2$-orthogonal projectors on these spaces are, respectively, denoted by $\RTproj{\ell}{Y}$ and $\NEproj{\ell}{Y}$.

\subsection{Discrete de Rham complex}

Throughout the rest of the paper, we fix an integer $k\ge 0$ corresponding to the polynomial degree of the DDR complex:
\begin{equation*}
  \begin{tikzcd}[column sep=2.5em]
    & 0 \arrow{r}{}
    & \Xgrad{k}{\Th}\arrow{r}{\uG{k}{h}}
    & \Xcurl{k}{\Th}\arrow{r}{\uC{k}{h}}
    & \Xdiv{k}{\Th}\arrow{r}{\D{k}{h}}
    &\Poly{k}(\Th)\arrow{r}{}
    & 0.
  \end{tikzcd}
\end{equation*}
In the following sections, we briefly recall the definitions of the spaces and the associated interpolators (which specify the meaning of the vector components in these spaces).
The components of the discrete operators connecting these spaces are defined by mimicking integration-by-parts formulae; we do not present these formulae as they are not directly useful to our purpose, and we refer the interested reader to \cite{Di-Pietro.Droniou:23} for details.

\subsubsection{Gradient space}

The discrete version of $\Hgrad{\Omega}$ is given in the DDR complex by
\[
\Xgrad{k}{\Th}\coloneq \bigtimes_{T\in\Th}\Poly{k-1}(T)\times\bigtimes_{F\in\Fh}\Poly{k-1}(F)\times \bigtimes_{E\in\Eh}\Poly{k-1}(E)\times \bigtimes_{V\in\Vh}\Real.
\]
The components of an element $\underline{q}_h\in\Xgrad{k}{\Th}$ are denoted by
\[
\underline{q}_h=((q_T)_{T\in\Th},(q_F)_{F\in\Fh},(q_E)_{E\in\Eh},(q_V)_{V\in\Vh}).
\]
The interpolator $\Igrad{k}{h}:C^0(\overline{\Omega})\to\Xgrad{k}{\Th}$ is defined by
\begin{equation*}
  \Igrad{k}{h}q
  \coloneq
  ((\lproj{k-1}{T}q)_{T\in\Th},(\lproj{k-1}{F}q)_{F\in\Fh},(\lproj{k-1}{E}q)_{E\in\Eh},(q(\bvec{x}_V))_{V\in\Vh})\quad\forall q\in C^0(\overline{\Omega}).
\end{equation*}

\subsubsection{Curl space}

The discrete $\Hcurl{\Omega}$ space is
\[
\Xcurl{k}{\Th}\coloneq \bigtimes_{T\in\Th}\RT{k}(T)\times \bigtimes_{F\in\Fh}\RT{k}(F)\times\bigtimes_{E\in\Eh}\Poly{k}(E)
\]
with components of $\uvec{v}_h\in\Xcurl{k}{\Th}$ written as
\[
\uvec{v}_h=((\bvec{v}_{\RTsymbol,T})_{T\in\Th},(\bvec{v}_{\RTsymbol,F})_{F\in\Fh},(v_E)_{E\in\Eh}).
\]
The interpolator $\Icurl{k}{h}:\bvec{C}^0(\overline{\Omega})\to \Xcurl{k}{\Th}$ is defined by
\begin{equation*}
  \Icurl{k}{h}\bvec{v}
  \coloneq
  ( (\RTproj{k}{T}\bvec{v})_{T\in\Th},(\RTproj{k}{F}\bvec{v}_{{\rm t},F})_{F\in\Fh},(\lproj{k}{E}(\bvec{v}\cdot\tangent_E))_{E\in\Eh})\quad\forall\bvec{v}\in
  \bvec{C}^0(\overline{\Omega}),
\end{equation*}
where $\bvec{v}_{{\rm t},F}=\normal_F\times(\bvec{v}\times\normal_F)$ denotes the tangential component of $\bvec{v}$ on $F$.

\subsubsection{Divergence space}

The counterpart of $\Hdiv{\Omega}$ in the DDR complex is
\[
\Xdiv{k}{\Th}\coloneq \bigtimes_{T\in\Th}\NE{k}(T)\times \bigtimes_{F\in\Fh}\Poly{k}(F)
\]
with a generic vector $\uvec{w}_h\in\Xdiv{k}{\Th}$ written as
\[
\uvec{w}_h=((\bvec{w}_{\NEsymbol,T})_{T\in\Th},(w_{F})_{F\in\Fh}).
\]
The interpolator $\Idiv{k}{h}:\bvec{C}^0(\overline{\Omega})\to\Xdiv{k}{\Th}$ is defined by
\begin{equation*}
  \Idiv{k}{h}\bvec{w}
  \coloneq
  ( (\NEproj{k}{T}\bvec{w})_{T\in\Th},(\lproj{k}{F}(\bvec{w}\cdot\normal_F))_{F\in\Fh} )\quad\forall\bvec{w}\in \bvec{C}^0(\overline{\Omega}).
\end{equation*}


\subsubsection{Discrete $L^2$-products and component $L^\alpha$-norms}

The restrictions of the DDR spaces to a mesh entity $Y\in\Eh\cup\Fh\cup\Th$ are denoted by replacing $\Th$ with $Y$, for example,
$\Xcurl{k}{F}\coloneq\RT{k}(F)\times\bigtimes_{E\in\EF}\Poly{k}(E)$,
while we replace the index $h$ by $Y$ to denote the restriction of vectors, as in $\uvec{v}_F= (\bvec{v}_{\RTsymbol,F},(v_E)_{E\in\EF}) \in \Xcurl{k}{F}$.

The Sobolev inequalities are stated using the following $L^\alpha$-like component norms on the discrete spaces: for $\alpha\in \lbrack 1,\infty\rparen$, $\bullet\in\{\GRAD,\CURL,\DIV\}$ and $\underline{x}_h=((x_T)_{T\in\Th},(x_F)_{F\in\Fh},(x_E)_{E\in\Eh},(x_V)_{V\in\Vh})\in\Xbullet{k}{\Th}$, with components taken equal to zero when they are not present in $\Xbullet{k}{\Th}$,
\begin{align*}
  \tnorm{\bullet,h,\alpha}{\underline{x}_h}&\coloneq \left(\sum_{T\in\Th}\tnorm{\bullet,T,\alpha}{\underline{x}_T}^\alpha\right)^{\frac1\alpha}\\
  \text{ where }\tnorm{\bullet,T,\alpha}{\underline{x}_T}&\coloneq\left(\norm{L^\alpha(T)}{x_T}^\alpha+h_T\sum_{F\in\FT}\norm{L^\alpha(F)}{x_F}^\alpha+h_T^2\sum_{E\in\ET}\norm{L^\alpha(E)}{x_E}^\alpha+h_T^3\sum_{V\in\VT}|x_V|^\alpha\right)^{\frac1\alpha}.
\end{align*}

Still considering $\bullet \in \{\GRAD,\CURL,\DIV\}$, set $k_{\bullet} \coloneqq k+1$ if $\bullet=\GRAD$, and $k_\bullet \coloneqq k$ if $\bullet\in\{\CURL,\DIV\}$. For all $T\in\Th$, a potential reconstruction is defined as $\Pbullet{k_\bullet}{T}:\Xbullet{k}{T}\to\Poly{k_\bullet}(T)$ (with a slight abuse of notation: for $\bullet\in\{\CURL,\DIV\}$ the codomain should be the vector-valued space $\vPoly{k}(T)$). These local potentials are patched together to form global piecewise discontinuous polynomial potentials:
\[
(\Pbullet{k_\bullet}{h}\underline{x}_h)_{|T}=\Pbullet{k_\bullet}{T}\underline{x}_T\quad\forall T\in\Th\,,\quad\forall\underline{x}_h\in\Xbullet{k}{\Th}.
\]
The $L^2$-like inner products on the DDR spaces are defined as follows: for all $\underline{x}_h, \underline{y}_h \in \Xbullet{k}{\Th}$,
\begin{equation*}
  \begin{aligned}
    (\underline{x}_h, \underline{y}_h)_{\bullet,h}
    &\coloneq \sum_{T \in \Th} (\underline{x}_T, \underline{y}_T)_{\bullet,T},
    \\
    \text{ where }
    (\underline{x}_T, \underline{y}_T)_{\bullet,T}
    &\coloneq \int_T \Pbullet{k_\bullet}{T}\underline{x}_T \cdot \Pbullet{k_\bullet}{T}\underline{y}_T
    + s_{\bullet,T}(\underline{x}_T,\underline{y}_T)\quad\forall T\in\Th,
  \end{aligned}
\end{equation*}
and the associated norms are denoted by $\norm{\bullet,T}{{\cdot}}$.
Here, the positive semidefinite stabilisation forms $s_{\bullet,T}$ must be chosen
(i) to scale in $h$ like the integral term in $(\cdot, \cdot)_{\bullet,T}$,
(ii) to vanish on interpolants of polynomials of degree $\le k_\bullet$ (polynomial consistency), and
(iii) to ensure the uniform local norm equivalence \cite[Lemma 5]{Di-Pietro.Droniou:23}:
\begin{equation}\label{eq:local.norm.equivalence}
  \norm{\bullet,T}{{\cdot}}\simeq \tnorm{\bullet,T,2}{{\cdot}}.
\end{equation}
These stabilisation forms are typically designed by penalising the difference between element potentials and potentials on faces and edges; see \cite[Section 4]{Di-Pietro.Droniou:23}.
The global stabilisations are obtained by summing the local ones:
\[
s_{\bullet,h}(\underline{x}_h,\underline{y}_h)=\sum_{T\in\Th}s_{\bullet,T}(\underline{x}_T,\underline{y}_T)\qquad
\forall \underline{x}_h,\underline{y}_h\in\Xbullet{k}{\Th},
\]
so that
\begin{equation}\label{eq:def.L2.product.global}
  (\underline{x}_h, \underline{y}_h)_{\bullet,h}=\int_\Omega \Pbullet{k_\bullet}{h}\underline{x}_h\cdot\Pbullet{k_\bullet}{h}\underline{y}_h + s_{\bullet,h}(\underline{x}_h,\underline{y}_h).
\end{equation}
The associated norm is $\norm{\bullet,h}{{\cdot}}$.

A straightforward consequence of \eqref{eq:local.norm.equivalence} is that the global
$L^2$-like norms and component norms are equivalent:
for all $\bullet\in\{\GRAD,\CURL,\DIV\}$,
\begin{equation}\label{eq:global.norm.equivalence}
  \norm{\bullet,h}{{\cdot}}\simeq \tnorm{\bullet,h,2}{{\cdot}}.
\end{equation}
Finally, the following local discrete Lebesgue inequality will also be of use (see \cite[Lemma 9]{Di-Pietro.Droniou.ea:24} for a proof): for all $\alpha,\beta\in \lbrack 1,\infty\rparen$ and $T\in\Th$,
\begin{equation}\label{eq:discrete.lebesgue}
  \tnorm{\bullet,T,\alpha}{\underline{x}_T}\lesssim h_T^{3\left(\frac1\alpha-\frac1\beta\right)}\tnorm{\bullet,T,\beta}{\underline{x}_T}\qquad\forall \underline{x}_T\in\Xbullet{k}{T}.
\end{equation}

\section{Sobolev inequalities}\label{sec:sobolev}

\subsection{Statement of the inequalities}

We state here the discrete Sobolev embeddings for DDR spaces, counterparts of \eqref{eq:continuous.sobolev.grad} and \eqref{eq:continuous.sobolev.curl.div}.

\begin{theorem}[Sobolev embedding in $\Xgrad{k}{\Th}$]\label{thm:sobolev.grad}
  The following inequality holds:
  \begin{equation}\label{eq:sobolev.grad}
    \tnorm{\GRAD,h,6}{\underline{q}_h}\lesssim \norm{\GRAD,h}{\underline{q}_h}+\norm{\CURL,h}{\uG{k}{h}\underline{q}_h}\qquad
    \forall \underline{q}_h\in \Xgrad{k}{\Th}.
  \end{equation}
\end{theorem}

  The embeddings for the curl and the divergence require the following mesh quasi-uniformity assumption:
  \begin{equation}\label{eq:quasi-uniform}
    h\lesssim h_T\qquad\forall T\in\Th.
  \end{equation}
  The hidden constants in the estimates and in the respective proofs may also depend on the quasi-uniformity constant.

  From this point on, the symbol $\perp$, used as a binary operator or a superscript, will denote orthogonality with respect to the appropriate (continuous or discrete) $L^2$-product.
  For example, in the following statement, $(\Image\uG{k}{h})^\perp$ denotes the orthogonal complement of $\Image\uG{k}{h}$ in $\Xcurl{k}{\Th}$ with respect to $(\cdot,\cdot)_{\CURL,h}$.

\begin{theorem}[Sobolev embedding in $\Xcurl{k}{\Th}$]\label{thm:sobolev.curl}
  Assume \eqref{eq:quasi-uniform}. Let $s\in (\frac12,1]$ satisfy \eqref{eq:regularity.Hcurl.Hdiv} and set $s^*=\frac{6}{3-2s}$. Then,
    \begin{equation}\label{eq:sobolev.curl}
      \tnorm{\CURL,h,s^*}{\uvec{v}_h}\lesssim  \norm{\CURL,h}{\uvec{v}_h}+\norm{\DIV,h}{\uC{k}{h}\uvec{v}_h}\qquad\forall \uvec{v}_h\in(\Image\uG{k}{h})^\perp.
    \end{equation}
\end{theorem}

\begin{theorem}[Sobolev embedding in $\Xdiv{k}{\Th}$]\label{thm:sobolev.div}
   Under the assumptions of Theorem \ref{thm:sobolev.curl}, the following inequality holds:
    \begin{equation}\label{eq:sobolev.div}
      \tnorm{\DIV,h,s^*}{\uvec{w}_h}\lesssim \norm{\DIV,h}{\uvec{w}_h}+\norm{L^2(\Omega)}{\D{k}{h}\uvec{w}_h}\qquad\forall \uvec{w}_h\in(\Image\uC{k}{h})^\perp
    \end{equation}
\end{theorem}

\begin{remark}[Bound on potential reconstructions]\label{rem:sobolev.potential}
  By \eqref{eq:local.norm.equivalence} together with the local polynomial Lebesgue inequalities \cite[Lemma 1.25]{Di-Pietro.Droniou:20} and \eqref{eq:discrete.lebesgue}, for $\bullet\in\{\GRAD,\CURL,\DIV\}$ and $\alpha\in [1,\infty)$ we have
    \[
    \norm{L^\alpha(\Omega)}{\Pbullet{k_\bullet}{h}\underline{x}_h}\lesssim \tnorm{\bullet,h,\alpha}{\underline{x}_h}\qquad\forall \underline{x}_h\in\Xbullet{k}{\Th}.
    \]
    The bounds in Theorems \ref{thm:sobolev.grad}--\ref{thm:sobolev.div} therefore directly give an $L^\alpha$-bound (with $\alpha=6$ or $\alpha=s^*$) on the potential reconstructions.
\end{remark}

\begin{remark}[Sobolev--Poincaré inequalities]
  Combining the theorems above with the DDR Poincaré inequalities \cite{Di-Pietro.Droniou:23,Di-Pietro.Hanot:24,Di-Pietro.Droniou.ea:25*1} gives the following discrete Sobolev--Poincaré bounds, in which the right-hand sides only depend on the vector through its discrete differential:
  \begin{alignat*}{2}
    \tnorm{\GRAD,h,6}{\underline{q}_h}&\lesssim  \norm{\CURL,h}{\uG{k}{h}\underline{q}_h}&&\quad\forall\underline{q}_h\in(\ker\uG{k}{h})^\perp,\\
    \tnorm{\CURL,h,s^*}{\uvec{v}_h}&\lesssim  \norm{\DIV,h}{\uC{k}{h}\uvec{v}_h}&&\quad\forall\uvec{v}_h\in(\ker\uC{k}{h})^\perp,\\
    \tnorm{\DIV,h,s^*}{\uvec{w}_h}&\lesssim  \norm{L^2(\Omega)}{\D{k}{h}\uvec{w}_h}&&\quad\forall\uvec{w}_h\in(\ker\D{k}{h})^\perp.
  \end{alignat*}
\end{remark}

\subsection{Proofs}

The proofs of the discrete Sobolev inequalities rely on the quasi-interpolators of \cite{Chaumont-Frelet.Droniou.ea:25} and the conforming lifting of \cite{Di-Pietro.Droniou.ea:25}. Theorem 7 in the former reference gives the existence, for each $\bullet\in\{\GRAD,\CURL,\DIV\}$, of a quasi-interpolator $\qIbullet{k}{h}:\Hbullet{\Omega}\to\Xbullet{k}{\Th}$ with the following locality and boundedness properties. For each $T\in\Th$, let $\widetilde{T}\coloneq \operatorname{int}\!\left(\bigcup\{\overline{T'}\,:\,T'\in\Th,\ \overline{T'}\cap\overline{T}\not=\emptyset\}\right)$ denote the patch of cells around $T$. The restriction $\qIbullet{k}{T}$ of $\qIbullet{k}{h}$ to $T$ has domain $\Hbullet{\widetilde{T}}$ and satisfies the following bound:
\begin{equation}
  \label{eq:qIbullet.bound}
  \tnorm{\bullet,T,2}{\qIbullet{k}{T} z}
  \lesssim \norm{{L}^2(\widetilde{T})}{z}+h_T\norm{\bvec{L}^2(\widetilde{T})}{\bullet z}\qquad\forall z\in\Hbullet{\widetilde{T}}.
\end{equation}
The quasi-interpolators commute with the differential operators so that, in particular,
\begin{equation}\label{eq:commutation.qI}
  \qIcurl{k}{h}(\GRAD q)= \uG{k}{h}(\qIgrad{k}{h} q)\qquad\forall q\in\Hgrad{\Omega}.
\end{equation}
By \cite[Lemma 18]{Chaumont-Frelet.Droniou.ea:25}, the quasi-interpolators can moreover be chosen in a way that is compatible with the lifting of \cite{Di-Pietro.Droniou.ea:25}, namely, such that each quasi-interpolator has a right inverse $\Liftbullet:\Xbullet{k}{\Th}\to\Hbullet{\Omega}$:
\begin{equation}
  \label{eq:qI.lift.inverse}
  \qIbullet{k}{h}(\Liftbullet\underline{z}_h)=\underline{z}_h\qquad\forall \underline{z}_h\in\Xbullet{k}{\Th},
\end{equation}
that satisfies the following bounds:
\begin{alignat}{2}
  \label{eq:bound.liftbullet}
  \norm{{L}^2(\Omega)}{\Liftbullet \underline{z}_h}&\lesssim \norm{\bullet,h}{\underline{z}_h}
  &&\qquad\forall\underline{z}_h\in\Xbullet{k}{\Th}\\
  \label{eq:bound.grad.liftgrad}
  \norm{\bvec{L}^2(\Omega)}{\GRAD\Liftgrad\underline{z}_h}&\lesssim \norm{\CURL,h}{\uG{k}{h}\underline{z}_h}
  &&\qquad\forall\underline{z}_h\in\Xgrad{k}{\Th},\\
  \label{eq:bound.curl.liftcurl}
  \norm{\bvec{L}^2(\Omega)}{\CURL\Liftcurl\uvec{v}_h}&\lesssim \norm{\DIV,h}{\uC{k}{h}\uvec{v}_h}
  &&\qquad\forall\uvec{v}_h\in\Xcurl{k}{\Th},\\
  \nonumber
  \norm{L^2(\Omega)}{\DIV\Liftdiv\uvec{w}_h}&\lesssim \norm{L^2(\Omega)}{\D{k}{h}\uvec{w}_h}
  &&\qquad\forall\uvec{w}_h\in\Xdiv{k}{\Th}.
\end{alignat}

\begin{proof}[Proof of Theorem \ref{thm:sobolev.grad} (Sobolev inequality for the gradient)]
  Let $q = \Liftgrad \underline{q}_h\in\Hgrad{\Omega}$.
  For all $T \in \Th$, we have
  \begin{align}
    \tnorm{\GRAD,T,6}{\qIgrad{k}{h} q}
    \overset{\eqref{eq:discrete.lebesgue}}&\lesssim h_T^{-1} \tnorm{\GRAD,T,2}{\qIgrad{k}{h} q} \nonumber\\
    \overset{\eqref{eq:qIbullet.bound}}&\lesssim h_T^{-1}\norm{L^2(\widetilde{T})}{q} + \norm{\bvec{L}^2(\widetilde{T})}{\GRAD q} \nonumber\\
    &\lesssim \norm{L^6(\widetilde{T})}{q} + \norm{\bvec{L}^2(\widetilde{T})}{\GRAD q}
    \label{eq:bound.L6.qI}
  \end{align}
  where the conclusion follows from a generalised H\"older inequality with exponents $(6,3)$ (which satisfy $\frac16+\frac13=\frac12$) together with the mesh regularity assumption, which ensures $|\widetilde{T}|^{\frac13}\lesssim h_T$.
  Note also that, for any family $(\xi_T)_{T\in\Th}$ of positive numbers and any $\zeta \ge 1$,
  \begin{equation}\label{eq:sum.xi.3}
    \sum_{T\in\Th}\xi_T^\zeta\le \left(\sum_{T\in\Th}\xi_T\right)^\zeta.
  \end{equation}
  We then write
  \begin{equation*}
    \begin{aligned}
      \tnorm{\GRAD,h,6}{\underline{q}_h}^6
      \overset{\eqref{eq:qI.lift.inverse}}&= \tnorm{\GRAD,h,6}{\qIgrad{k}{h} q}^6 \\
      \overset{\eqref{eq:bound.L6.qI}}&\lesssim \sum_{T\in\Th} \norm{L^6(\widetilde{T})}{q}^6
      + \sum_{T\in\Th} \left( \norm{\bvec{L}^2(\widetilde{T})}{\GRAD q}^2 \right)^3 \\
      \overset{\eqref{eq:sum.xi.3}}&\lesssim \norm{L^6(\Omega)}{q}^6
      + \left(\sum_{T\in\Th}  \norm{\bvec{L}^2(\widetilde{T})}{\GRAD q}^2 \right)^3\\
      &\lesssim \norm{L^6(\Omega)}{q}^6
      + \norm{\bvec{L}^2(\Omega)}{\GRAD q}^6,
    \end{aligned}
  \end{equation*}
  where we have used, in the second and third inequalities, the fact that each $T'\in\Th$ appears in $\lesssim 1$ patches $\widetilde{T}$ to write, for all $\rho\in [1,\infty)$ and $z\in L^\rho(\Omega)$,
    \begin{equation}\label{eq:norm.patches.Lrho}
      \sum_{T\in\Th}\norm{{L}^\rho(\widetilde{T})}{z}^\rho\lesssim\sum_{T\in\Th}\norm{L^\rho(T)}{z}^\rho=\norm{L^\rho(\Omega)}{z}^\rho.
    \end{equation}
    Invoking the Sobolev embedding $\Hgrad{\Omega}\hookrightarrow L^6(\Omega)$ then yields
    \[
    \begin{aligned}
      \tnorm{\GRAD,h,6}{\underline{q}_h}^6 &\lesssim \norm{L^2(\Omega)}{q}^6 + \norm{\bvec{L}^2(\Omega)}{\GRAD q}^6\\
      \overset{\eqref{eq:bound.liftbullet}, \eqref{eq:bound.grad.liftgrad}}&\lesssim \norm{\GRAD,h}{\underline{q}_h}^{6} + \norm{\CURL,h}{\uG{k}{h}\underline{q}_h}^{6},
    \end{aligned}
    \]
    which concludes the proof.
\end{proof}

\begin{proof}[Proof of Theorem \ref{thm:sobolev.curl} (Sobolev inequality for the curl)]
  Let $\uvec{v}_h\in(\Image\uG{k}{h})^\perp$. The Hodge decomposition of $\Liftcurl\uvec{v}_h\in\Hcurl{\Omega}$ along with the $L^2$-orthogonal decomposition $\Hcurl{\Omega}=(\Hcurl{\Omega}\cap(\Image\GRAD)^\perp)\oplus\Image\GRAD$ gives
  \begin{equation}\label{eq:Lift.hodge}
    \text{$\Liftcurl\uvec{v}_h=\bvec{w}+\GRAD q$ with $\bvec{w}\in\Hcurl{\Omega}\cap\Hdivzz{\Omega}$ and $q\in\Hgrad{\Omega}$},
  \end{equation}
  where we have used the fact that
  \begin{equation}\label{eq:im.grad.div.0}
    (\Image\GRAD)^\perp=\Hdivzz{\Omega}\coloneq\{\bvec{z}\in\Hdiv{\Omega}\,:\,\DIV\bvec{z}=0\,,\;\bvec{z}\cdot\normal_\Omega=0\text{ on $\partial\Omega$}\}.
  \end{equation}
  Since the decomposition in \eqref{eq:Lift.hodge} is $L^2$-orthogonal and $\CURL\Liftcurl\uvec{v}_h=\CURL\bvec{w}$, we have
  \begin{equation}\label{eq:bound.w.q}
    \begin{aligned}
      \norm{\bvec{L}^2(\Omega)}{\bvec{w}}+\norm{\bvec{L}^2(\Omega)}{\GRAD q}\le 2\norm{\bvec{L}^2(\Omega)}{\Liftcurl\uvec{v}_h}
      \overset{\eqref{eq:bound.liftbullet}}&\lesssim \norm{\CURL,h}{\uvec{v}_h},\\
      \norm{\bvec{L}^2(\Omega)}{\CURL\bvec{w}}= \norm{\bvec{L}^2(\Omega)}{\CURL\Liftcurl\uvec{v}_h}\overset{\eqref{eq:bound.curl.liftcurl}}&\lesssim\norm{\DIV,h}{\uC{k}{h}\uvec{v}_h}.
    \end{aligned}
  \end{equation}
  Moreover, \eqref{eq:regularity.Hcurl.Hdiv} gives $\bvec{w}\in \bvec{H}^s(\Omega)$. Since $\DIV\bvec{w}=0$,
  \begin{equation}\label{eq:w.H1}
    \norm{\bvec{H}^s(\Omega)}{\bvec{w}}\lesssim \norm{\bvec{L}^2(\Omega)}{\bvec{w}}+\norm{\bvec{L}^2(\Omega)}{\CURL\bvec{w}}\overset{\eqref{eq:bound.w.q}}\lesssim  \norm{\CURL,h}{\uvec{v}_h}+\norm{\DIV,h}{\uC{k}{h}\uvec{v}_h}.
  \end{equation}

  Applying the quasi-interpolator $\qIcurl{k}{h}$ to the decomposition \eqref{eq:Lift.hodge}, using the right-inverse property \eqref{eq:qI.lift.inverse} and invoking the commutation property \eqref{eq:commutation.qI}, we obtain
  \begin{equation}\label{eq:decomp.vh}
    \uvec{v}_h=\qIcurl{k}{h}\bvec{w}+\uG{k}{h}\qIgrad{k}{h}q.
  \end{equation}
  The estimate \eqref{eq:sobolev.curl} is obtained by bounding each term on the right-hand side above. Let us start with the gradient term.
  Since, by assumption, $\uvec{v}_h$ is $(\cdot,\cdot)_{\CURL,h}$-orthogonal to $\uG{k}{h}\qIgrad{k}{h}q \in \Image \uG{k}{h} $, we have
  \[
  0=(\uvec{v}_h,\uG{k}{h}\qIgrad{k}{h}q)_{\CURL,h}\overset{\eqref{eq:decomp.vh}}=(\qIcurl{k}{h}\bvec{w},\uG{k}{h}\qIgrad{k}{h}q)_{\CURL,h}+\norm{\CURL,h}{\uG{k}{h}\qIgrad{k}{h}q}^2.
  \]
  Since $\bvec{w}\in\Hdiv{\Omega}$ with $\DIV\bvec{w}=0$ and $\bvec{w}\cdot\normal_\Omega=0$ on $\partial\Omega$, the adjoint consistency result of \cite[Theorem 11, (i)]{Chaumont-Frelet.Droniou.ea:25} allows us to continue with
  \begin{equation}
    \begin{aligned}
      &\norm{\CURL,h}{\uG{k}{h}\qIgrad{k}{h}q}^2
      \\
      &\quad
      \begin{aligned}[t]
        &=-(\qIcurl{k}{h}\bvec{w},\uG{k}{h}\qIgrad{k}{h}q)_{\CURL,h}
        \\
        &\lesssim \left(\norm{\bvec{L}^2(\Omega)}{\vlproj{k}{h}\bvec{w}-\bvec{w}}+h\norm{\bvec{L}^2(\Omega)}{\vlproj{k}{h}(\CURL\bvec{w})-\CURL\bvec{w}}\right)\tnorm{\CURL,h,2}{\uG{k}{h}\qIgrad{k}{h}q}\\
        &\lesssim h^s \left(\norm{\CURL,h}{\uvec{v}_h}+\norm{\DIV,h}{\uC{k}{h}\uvec{v}_h}\right)\tnorm{\CURL,h,2}{\uG{k}{h}\qIgrad{k}{h}q},
      \end{aligned}
    \end{aligned}
    \label{eq:use.adjoint}
  \end{equation}
  where we have used the approximation and boundedness properties of local $L^2$-projectors on polynomial spaces together with \eqref{eq:w.H1}.
  Using the norm equivalence \eqref{eq:global.norm.equivalence} and simplifying leads to
  \begin{equation}\label{eq:est.Ghk.q.L2}
    \tnorm{\CURL,h,2}{\uG{k}{h}\qIgrad{k}{h}q}^2
    \lesssim h^{2s} \left(\norm{\CURL,h}{\uvec{v}_h}+\norm{\DIV,h}{\uC{k}{h}\uvec{v}_h}\right)^2.
  \end{equation}
  We then invoke the discrete Lebesgue inequality \eqref{eq:discrete.lebesgue} with $(\alpha,\beta)=(s^*,2)$ and observe that $3\left(\frac{1}{s^*}-\frac12\right)=\frac{3-2s-3}{2}=-s$ to write
  \begin{align}
    \tnorm{\CURL,h,s^*}{\uG{k}{h}\qIgrad{k}{h}q}^{s^*}&=\sum_{T\in\Th}\tnorm{\CURL,T,s^*}{\uG{k}{T}\qIgrad{k}{T}q}^{s^*}\nonumber\\
    &\lesssim \sum_{T\in\Th}h_T^{-s^*s}\tnorm{\CURL,T,2}{\uG{k}{T}\qIgrad{k}{T}q}^{s^*}\nonumber\\
    \overset{\eqref{eq:quasi-uniform}}&\lesssim h^{-s^*s}\left(\sum_{T\in\Th}\tnorm{\CURL,T,2}{\uG{k}{T}\qIgrad{k}{T}q}^2\right)^{\frac{s^*}{2}}\nonumber\\
    \overset{\eqref{eq:est.Ghk.q.L2}}&\lesssim \left(\norm{\CURL,h}{\uvec{v}_h}+\norm{\DIV,h}{\uC{k}{h}\uvec{v}_h}\right)^{s^*}
    \label{eq:est.Ghk.q.L6}
  \end{align}
  where, to pass to the second inequality, we have additionally used \eqref{eq:sum.xi.3} with $\zeta = \frac{s^*}{2} \ge 1$.

  We now turn to the first term on the right-hand side of \eqref{eq:decomp.vh}. Let $T\in\Th$. Using again the discrete Lebesgue inequality \eqref{eq:discrete.lebesgue} with $(\alpha,\beta)=(s^*,2)$ followed by the bound \eqref{eq:qIbullet.bound} on $\qIcurl{k}{T}$, we have
  \begin{align*}
    \tnorm{\CURL,T,s^*}{\qIcurl{k}{T}\bvec{w}}&\lesssim h_{T}^{-s}\tnorm{\CURL,T,2}{\qIcurl{k}{T}\bvec{w}}\nonumber\\
    &\lesssim h_{T}^{-s}\norm{\bvec{L}^2(\widetilde{T})}{\bvec{w}}+h_T^{1-s}\norm{\bvec{L}^2(\widetilde{T})}{\CURL\bvec{w}}\nonumber\\
    &\lesssim \norm{\bvec{L}^{s^*}(\widetilde{T})}{\bvec{w}}+\norm{\bvec{L}^2(\widetilde{T})}{\CURL\bvec{w}}
  \end{align*}
  where the conclusion follows from $h_T^{1-s}\le \mathrm{diam}(\Omega)^{1-s}$ (since $1-s\ge 0$) and, as in \eqref{eq:bound.L6.qI}, a generalised H\"older inequality (with exponents $\big(s^*,\frac{3}{s}\big)$, which satisfy $\frac{1}{s^*}+\frac{s}{3}=\frac12$) and the mesh regularity assumption. We raise this estimate to the power $s^*$, use $(a+b)^{s^*}\lesssim a^{s^*}+b^{s^*}$ on the right-hand side, sum over $T\in\Th$
  and invoke \eqref{eq:sum.xi.3} with $\xi_T=\norm{\bvec{L}^2(\widetilde{T})}{\CURL\bvec{w}}^2$ and $\zeta=\frac{s^*}{2}$ to get
  \begin{align}
    \tnorm{\CURL,h,s^*}{\qIcurl{k}{h}\bvec{w}}^{s^*}&\lesssim \sum_{T\in\Th}\norm{\bvec{L}^{s^*}(\widetilde{T})}{\bvec{w}}^{s^*}
    +\left(\sum_{T\in\Th}\norm{\bvec{L}^2(\widetilde{T})}{\CURL\bvec{w}}^2\right)^{\frac{s^*}{2}}\nonumber\\
    \overset{\eqref{eq:norm.patches.Lrho}}&\lesssim \norm{\bvec{L}^{s^*}(\Omega)}{\bvec{w}}^{s^*}+\norm{\bvec{L}^2(\Omega)}{\CURL\bvec{w}}^{s^*}\nonumber\\
    \overset{\eqref{eq:w.H1}}&\lesssim \left(\norm{\CURL,h}{\uvec{v}_h}+\norm{\DIV,h}{\uC{k}{h}\uvec{v}_h}\right)^{s^*},
    \label{eq:est.qIw}
  \end{align}
  where we have additionally used the Sobolev embedding $H^s(\Omega)\hookrightarrow L^{s^*}(\Omega)$ in the conclusion.

  The inequality \eqref{eq:sobolev.curl} is obtained by applying $\tnorm{\CURL,h,s^*}{{\cdot}}$ to \eqref{eq:decomp.vh}, using the triangle inequality, and plugging in \eqref{eq:est.Ghk.q.L6} and \eqref{eq:est.qIw} raised to the power $\frac{1}{s^*}$.
\end{proof}

\begin{proof}[Proof of Theorem \ref{thm:sobolev.div} (Sobolev inequality for the divergence)]
  The proof follows the same argument as that of Theorem \ref{thm:sobolev.curl}, using the lifting and quasi-interpolators in $\Xdiv{k}{\Th}$ and relying on the $L^2$-orthogonal Hodge decomposition
  \begin{align*}
    \Hdiv{\Omega}&=(\Hdiv{\Omega}\cap(\Image\CURL)^\perp)\oplus\Image\CURL\\
    &=\{\bvec{z}\in\Hdiv{\Omega}\cap\Hcurl{\Omega}\,:\,\CURL\bvec{z}=0\text{ and }\bvec{z}\times\normal_\Omega=0\text{ on $\partial\Omega$}\}\oplus\Image\CURL
    \qedhere
  \end{align*}
\end{proof}


\section{Application to the analysis of the DDR method for the Navier--Stokes equations}\label{sec:ddr.ns}

\subsection{Model and scheme}

On a convex domain, the embedding \eqref{eq:regularity.Hcurl.Hdiv} is valid for $s=1$, and the identities
\[
-\Delta \bvec{u}=\CURL\CURL\bvec{u}-\GRAD\DIV\bvec{u}\text{ and } \DIV(\bvec{u}\otimes\bvec{u})=(\DIV\bvec{u})\bvec{u}+(\CURL \bvec{u})\times\bvec{u}+\frac12\GRAD(\bvec{u}\cdot\bvec{u})
\]
therefore formally justify the $\CURL\CURL$ formulation of the Navier--Stokes equations \cite{Girault:90} which, for $\bvec{f}\in \bvec{L}^{\frac{4}{3}}(\Omega)$, reads:
Find $(\bvec{u},p)\in\Hcurl{\Omega}\times W^{1,\frac{4}{3}}(\Omega)$ such that $\int_\Omega p=0$ and
\begin{subequations}\label{eq:ns.curl}
  \begin{alignat}{2}
    \label{eq:ns.curl.momentum}
    &\nu\int_\Omega\CURL\bvec{u}\cdot\CURL\bvec{v}+\int_\Omega [(\CURL\bvec{u})\times\bvec{u}]\cdot\bvec{v}+\int_\Omega\GRAD p\cdot\bvec{v}=\int_\Omega \bvec{f}\cdot\bvec{v}&&\quad\forall \bvec{v}\in\Hcurl{\Omega}\cap\bvec{L}^4(\Omega),\\
    \label{eq:ns.curl.incompressibility}
    &-\int_\Omega \bvec{u}\cdot\GRAD q=0&&\quad\forall q\in\Hgrad{\Omega},
  \end{alignat}
\end{subequations}
where $\nu>0$ is the viscosity and $\bvec{f}$ the forcing term. These equations correspond to the following (natural) boundary
conditions:
\[
\text{$\CURL\bvec{u}\times\bvec{n}=0$ and $\bvec{u}\cdot\bvec{n}=0$ on $\partial\Omega$}.
\]
For simplicity, we use these boundary conditions below, but we note that, as discussed in \cite[Section 5]{Di-Pietro.Droniou.ea:24}, the essential boundary conditions
\[
\text{$\bvec{u}\times\bvec{n}=0$ and $p=0$ on $\partial\Omega$}
\]
could also be considered by using $\Hcurlz{\Omega}$ and $\Hgradz{\Omega}$ as test and trial spaces in \eqref{eq:ns.curl}.
Assuming that $\bvec{f}\in \bvec{C}(\overline{\Omega})$, the DDR scheme for \eqref{eq:ns.curl} reads: Find $(\uvec{u}_h,\underline{p}_h)\in\Xcurl{k}{\Th}\times\Xgrad{k}{\Th}$ such that $(\underline{p}_h,\Igrad{k}{h}1)_{\GRAD,h}=0$ and
\begin{subequations}\label{eq:ddr.scheme}
  \begin{alignat}{2}
    \label{eq:ddr.scheme.momentum}
    &\begin{aligned}
       \nu(\uC{k}{h}\uvec{u}_h,\uC{k}{h}\uvec{v}_h)_{\DIV,h}
       +\int_\Omega {}&[(\Pdiv{k}{h}\uC{k}{h}\uvec{u}_h)\times\Pcurl{k}{h}\uvec{u}_h]\cdot\Pcurl{k}{h}\uvec{v}_h\\
       &+(\uG{k}{h} \underline{p}_h,\uvec{v}_h)_{\CURL,h}=(\Icurl{k}{h}\bvec{f},\uvec{v}_h)_{\CURL,h}
     \end{aligned}&&\qquad\forall \uvec{v}_h\in\Xcurl{k}{\Th},\\
    \label{eq:ddr.scheme.incompressibility}
    &-(\uvec{u}_h,\uG{k}{h}\underline{q}_h)_{\CURL,h}=0&&\qquad\forall \underline{q}_h\in\Xgrad{k}{\Th}.
  \end{alignat}
\end{subequations}

  \begin{remark}[Comparison with the scheme of \cite{Di-Pietro.Droniou.ea:24}]
    In contrast to the scheme in \cite{Di-Pietro.Droniou.ea:24}, here we do not use serendipity to reduce the number of element and face unknowns.
    This corresponds to taking $\eta_Y = 2$ for all $Y \in \Th \cup \Fh$ in the notation of this reference.
    Both the serendipity and non-serendipity DDR complexes share the same relevant properties \cite{Di-Pietro.Droniou:23*1}; hence the analysis that follows can be extended in a straightforward way to cover the serendipity case.
  \end{remark}

  \begin{proposition}[Existence of a discrete solution]\label{prop:existence.discrete}
    If $\Omega$ is convex, $\nu>0$ and $\bvec{f}\in\bvec{C}(\overline{\Omega})$, then \eqref{eq:ddr.scheme} has at least one solution on every admissible mesh.
  \end{proposition}
  
  \begin{proof}
    Set $\cvec{Z}_h \coloneqq (\Image\uG{k}{h})^\perp$.
    Exactness on the convex domain and the discrete Poincar\'e inequality \cite{Di-Pietro.Droniou.ea:25*1} give
    $\norm{\CURL,h}{\uvec{z}_h}\le \Cpoin \norm{\DIV,h}{\uC{k}{h}\uvec{z}_h}$ for all $\uvec{z}_h\in\cvec{Z}_h$.
    Define $A_h:\cvec{Z}_h\to\cvec{Z}_h$ such that, for all $\uvec{z}_h,\uvec{v}_h\in\cvec{Z}_h$,
    \[
      \begin{aligned}
        (A_h\uvec{z}_h,\uvec{v}_h)_{\CURL,h}
        &=\nu(\uC{k}{h}\uvec{z}_h,\uC{k}{h}\uvec{v}_h)_{\DIV,h}\\
        &\quad
        +\int_\Omega[(\Pdiv{k}{h}\uC{k}{h}\uvec{z}_h)\times\Pcurl{k}{h}\uvec{z}_h]\cdot\Pcurl{k}{h}\uvec{v}_h
        -(\Icurl{k}{h}\bvec{f},\uvec{v}_h)_{\CURL,h}.
      \end{aligned}
    \]
    The right-hand side is linear in $\uvec{v}_h$, so the Riesz representation theorem uniquely defines $A_h\uvec{z}_h$.
    Its dependence on $\uvec{z}_h$ is polynomial, so $A_h$ is continuous.
    Selecting $\uvec{v}_h = \uvec{z}_h$ and using the skew-symmetry of the nonlinear term, the Poincar\'e bound on $\uvec{z}_h$ and the Cauchy--Schwarz inequality for the source term yields
    \[
    (A_h \uvec{z}_h,\uvec{z}_h)_{\CURL,h}
    \ge \nu\Cpoin^{-2}\norm{\CURL,h}{\uvec{z}_h}^{2}
    -\norm{\CURL,h}{\Icurl{k}{h}\bvec{f}}\norm{\CURL,h}{\uvec{z}_h}.
    \]
    This estimate shows that $(A_h \uvec{z}_h,\uvec{z}_h)_{\CURL,h}/\norm{\CURL,h}{\uvec{z}_h}\to \infty$ as $\norm{\CURL,h}{\uvec{z}_h}\to\infty$ and thus, by \cite[Theorem 3.3]{Deimling:85}, that $A_h$ has a zero $\uvec{u}_h\in\mathcal Z_h$.
    \smallskip
    
    For this velocity, define the remaining momentum residual $\mathcal{R}_h:\Xcurl{k}{\Th}\to\Real$ by
    \[
    \mathcal{R}_h(\uvec{v}_h)\coloneqq
    (\Icurl{k}{h}\bvec{f},\uvec{v}_h)_{\CURL,h}
    -\nu(\uC{k}{h}\uvec{u}_h,\uC{k}{h}\uvec{v}_h)_{\DIV,h}
    -\int_\Omega[(\Pdiv{k}{h}\uC{k}{h}\uvec{u}_h)\times\Pcurl{k}{h}\uvec{u}_h]\cdot\Pcurl{k}{h}\uvec{v}_h.
    \]
    This is a linear functional on the full velocity space. For every $\uvec{v}_h\in\cvec{Z}_h$, the definition of $A_h$ gives
    $\mathcal{R}_h(\uvec{v}_h)=-(A_h\uvec{u}_h,\uvec{v}_h)_{\CURL,h}=0$.
    Its Riesz representative therefore belongs to $\cvec{Z}_h^\perp=\Image\uG{k}{h}$, so there exists $\underline{p}_h\in\Xgrad{k}{\Th}$ such that
    \[
      (\uG{k}{h}\underline{p}_h,\uvec{v}_h)_{\CURL,h}
      =\mathcal{R}_h(\uvec{v}_h)
      \qquad\forall\uvec{v}_h\in\Xcurl{k}{\Th}.
    \]
    This is exactly \eqref{eq:ddr.scheme.momentum}, while $\uvec{u}_h\in\cvec{Z}_h$ gives \eqref{eq:ddr.scheme.incompressibility}.
    Since $\ker\uG{k}{h}= \Igrad{k}{h} \Real$, the pressure is uniquely determined by the additional condition $(\underline{p}_h,\Igrad{k}{h}1)_{\GRAD,h}=0$.
  \end{proof}

\subsection{Error estimate for smooth and small data}

In \cite{Di-Pietro.Droniou.ea:24}, the scheme \eqref{eq:ddr.scheme} was analysed under the \emph{assumption} that a uniform-in-$h$ Sobolev inequality holds in $\Xcurl{k}{\Th}$, and using the following discrete seminorm, which is a norm on the normalised pressure space, corresponding to a dual $L^{\frac{4}{3}}$-norm on the gradients:
\begin{equation}\label{eq:def.norm.P}
\norm{P}{\underline{r}_h}\coloneq \max\biggl\{(\uG{k}{h}\underline{r}_h,\uG{k}{h}\underline{q}_h)_{\CURL,h}\,:
    \, \underline{q}_h\in\Xgrad{k}{\Th}\,,\;\tnorm{\CURL,h,4}{\uG{k}{h}\underline{q}_h}\le 1\biggr\}.
\end{equation}

Theorem \ref{thm:sobolev.curl} establishes the required uniform Sobolev inequality and allows us to sharpen the error estimate from \cite[Theorem 2]{Di-Pietro.Droniou.ea:24} as follows:

\begin{theorem}[Error estimate]\label{thm:error.estimates}
  Assume that $\Omega$ is convex and that the meshes satisfy \eqref{eq:quasi-uniform}.
  There exists $\eta>0$ depending only on $\Omega$, the mesh regularity parameter, the polynomial degree $k$, and the mesh quasi-uniformity constant in \eqref{eq:quasi-uniform}, such that, if \eqref{eq:ns.curl} has a solution $(\bvec{u},p)\in\Hcurl{\Omega}\times\Hgrad{\Omega}$ satisfying
  \begin{equation}\label{eq:reg.sol.error.estimate}
  \begin{aligned}
    &\bvec{u}\in \bvec{W}^{1,4}(\Omega)\cap \bvec{W}^{k+1,4}(\Th)\cap \bvec{H}^{\max(k+1,2)}(\Th),\\
    &\CURL\bvec{u}\in \bvec{C}^0(\overline{\Omega})\cap \bvec{H}^{k+2}(\Th)\cap \bvec{W}^{k+1,4}(\Th),\\
    &(\CURL \bvec{u})\times\bvec{u}\in \bvec{H}^{\max(k+1,2)}(\Th)\,,\quad\CURL\CURL\bvec{u}\in \bvec{H}^{\max(k+1,2)}(\Th)
  \end{aligned}
  \end{equation}
  and
  \begin{equation}\label{eq:smallness}
    \chi\coloneq \nu-\eta\nu^{-1}\norm{\CURL,h}{\Icurl{k}{h}[\nu\CURL\CURL\bvec{u} + (\CURL\bvec{u})\times\bvec{u}]}>0,
  \end{equation}
  then any solution $(\uvec{u}_h,\underline{p}_h)\in\Xcurl{k}{\Th}\times \Xgrad{k}{\Th}$ to \eqref{eq:ddr.scheme} satisfies
  \begin{equation}\label{eq:error.u}
    \norm{\CURL,h}{\uvec{u}_h-\Icurl{k}{h}\bvec{u}}
    +\norm{\DIV,h}{\uC{k}{h}(\uvec{u}_h-\Icurl{k}{h}\bvec{u})}
    +\norm{P}{\underline{p}_h-\Igrad{k}{h}p}
    \le \chi^{-1} \mathcal K(\bvec{u}) h^{k+1},
\end{equation}
  where $\mathcal K(\bvec{u})$ only depends on $\bvec{u}$.
\end{theorem}

\begin{remark}[Pressure-robust error estimates]
  As shown in \cite{Di-Pietro.Droniou.ea:24}, the constant $\eta$ in Theorem \ref{thm:error.estimates} can be explicitly evaluated in terms of continuity constants of the DDR potentials, discrete Poincaré constants and the Sobolev inequality constant with exponent $4$.
 
  Note also that \eqref{eq:error.u} is robust with respect to the pressure: the constant on the right-hand side does not depend on $p$,
  and the smallness assumption \eqref{eq:smallness} only involves the velocity-dependent part of the forcing term $\bvec{f}=\nu\CURL\CURL\bvec{u}+(\CURL\bvec{u})\times\bvec{u}+\GRAD p$.
\end{remark}

\begin{remark}[Regularity of $\bvec{f}$]
  The continuity of $\bvec{f}$ is assumed in order to use its interpolant to define the right-hand side in \eqref{eq:ddr.scheme.momentum}. If we only suppose that $\bvec{f}\in \bvec{L}^{\frac{4}{3}}(\Omega)$, then a scheme can be written by replacing this right-hand side with $\int_\Omega \bvec{f}\cdot\Pcurl{k}{h}\uvec{v}_h$. Under the regularity assumptions \eqref{eq:reg.sol.error.estimate}, this scheme would still have optimal $\mathcal O(h^{k+1})$ convergence rates; it would however not be pressure-robust.

  The convergence analysis by compactness presented in Section \ref{sec:convergence.compactness} would also apply (without smoothness assumptions on the data) with this modified right-hand side.
\end{remark}

\subsection{Convergence by compactness for general data and solutions}\label{sec:convergence.compactness}

The error estimate in Theorem \ref{thm:error.estimates} for the stationary Navier--Stokes equations requires both a smallness assumption on the data and smoothness of the velocity. Convergence of the scheme can nevertheless be proved without these assumptions, provided that sufficient compactness is available for the discrete solutions. Here, this is ensured by the recent Maxwell compactness result for DDR proved in \cite{Chaumont-Frelet.Droniou.ea:25}.

In estimates for the discrete solutions in the following compactness analysis, hidden constants may additionally depend on the fixed viscosity $\nu>0$ and forcing term $\bvec{f}$.

\begin{theorem}[Convergence by compactness]\label{thm:convergence.compactness}
  Assume that $\Omega$ is convex and that $\bvec{f}\in\bvec{C}(\overline{\Omega})$. Let $(\Mh)_{h\in\mathcal H}$ be a sequence of quasi-uniform regular meshes and, for all $h\in\mathcal H$, let $(\uvec{u}_h,\underline{p}_h)$ be a solution to \eqref{eq:ddr.scheme}.
  Then, there exists a solution $(\bvec{u},p)\in\Hcurl{\Omega}\times  W^{1,\frac{4}{3}}(\Omega)$ to \eqref{eq:ns.curl} such that, as $h\to 0$ and along a subsequence (not relabelled), we have
  \begin{alignat}{2}
    \label{eq:conv.Pcurlh.u}
    \Pcurl{k}{h}\uvec{u}_h&\to \bvec{u}&&\quad\text{strongly in $\bvec{L}^{r}(\Omega)$ for all $r \in [1,6)$},\\
      \label{eq:conv.Ch.u}
      \Pdiv{k}{h}\uC{k}{h}\uvec{u}_h&\to \CURL\bvec{u}&&\quad\text{strongly in $\bvec{L}^{2}(\Omega)$},\\
      \label{eq:conv.stab.Ch.u}
      s_{\DIV,h}(\uC{k}{h}\uvec{u}_h,\uC{k}{h}\uvec{u}_h)&\to 0,&&\\
      \label{eq:conv.Pgradh.p}
      \Pgrad{k+1}{h}\underline{p}_h&\to p&&\quad\text{strongly in $L^{2}(\Omega)$}.
  \end{alignat}
\end{theorem}

\begin{proof} See Section \ref{sec:proof.convergence.compactness}.
\end{proof}

\subsubsection{Preliminary results}

The proof of Theorem \ref{thm:convergence.compactness} relies on several technical results that we establish here. The first one is a compactness lemma in $\Xcurl{k}{\Th}$, the second states an approximation property for functions in $\Hcurl{\Omega}\cap(\Image\GRAD)^\perp$, the third concerns bounds derived from the discrete Lebesgue inequality \eqref{eq:discrete.lebesgue}, and the last concerns the properties of the interpolator and quasi-interpolator in non-Hilbertian settings.

\begin{lemma}[Compactness in $\Xcurl{k}{\Th}$ and regularity of the limit]\label{lem:compactness.curl}
  Assume that $\Omega$ is convex and let $(\Mh)_{h\in\mathcal H}$ be a sequence of quasi-uniform regular meshes. For all $h\in\mathcal H$, take $\uvec{u}_h\in\Xcurl{k}{\Th}$ such that $\uvec{u}_h\perp \ker\uC{k}{h}$ and $(\norm{\DIV,h}{\uC{k}{h}\uvec{u}_h})_{h\in\mathcal H}$ is bounded. Then, there exists $\bvec{u}\in \Hcurl{\Omega}$ such that $\bvec{u}\perp \ker\CURL$ and, up to a subsequence as $h\to 0$,
  \begin{alignat}{2}
  \label{eq:conv.Ph.u.strong}
   \Pcurl{k}{h}\uvec{u}_h&\to \bvec{u}&&\quad\text{strongly in $\bvec{L}^{r}(\Omega)$ for all $r \in [1,6)$},\\
  \label{eq:conv.Ch.u.weak}
    \Pdiv{k}{h}\uC{k}{h}\uvec{u}_h&\to \CURL\bvec{u}&&\quad\text{weakly in $\bvec{L}^2(\Omega)$.}
  \end{alignat}
\end{lemma}

\begin{remark}[Convexity and quasi-uniformity]
  The convexity of $\Omega$ and the quasi-uniformity of $\Mh$ are only used to invoke the discrete Sobolev inequalities and obtain the convergence \eqref{eq:conv.Pcurlh.u} for all $r\in [1,6)$. Without these assumptions, Lemma \ref{lem:compactness.curl} still holds with $r=2$.
\end{remark}

\begin{proof}
  The discrete Poincar\'e inequality \cite{Di-Pietro.Droniou.ea:25*1} and the orthogonality assumption yield
  \begin{equation}
    \norm{\CURL,h}{\uvec{u}_h}
    \lesssim \norm{\DIV,h}{\uC{k}{h}\uvec{u}_h}
    \lesssim 1.
    \label{eq:conv.Ch.u.weak:basic-bounds}
  \end{equation}
  Consequently, by \eqref{eq:global.norm.equivalence} and Remark \ref{rem:sobolev.potential}, $(\Pdiv{k}{h}\uC{k}{h}\uvec{u}_h)_{h\in\mathcal H}$ is bounded in
  $\bvec{L}^2(\Omega)$.
  Hence, up to a subsequence, there exists $\bvec{\xi}\in\bvec{L}^2(\Omega)$ such that
  \begin{equation}\label{eq:conv.Ch.u.xi}
    \Pdiv{k}{h}\uC{k}{h}\uvec{u}_h
    \to\bvec{\xi} \quad\text{weakly in $\bvec{L}^2(\Omega)$}.
  \end{equation}
  Moreover, applying \cite[Theorem 2 and Remark 5]{Chaumont-Frelet.Droniou.ea:25} in the unweighted case and with no boundary conditions, we infer the
  existence of $\bvec{u}\in\bvec{L}^2(\Omega)$ such that, up to a subsequence as $h\to 0$,
  \begin{gather}\label{eq:conv.Pcurlh.u.L2}
    \Pcurl{k}{h}\uvec{u}_h
    \to\bvec{u} \quad\text{ strongly in $\bvec{L}^2(\Omega)$},\\
    \label{eq:conv.stab.curlh.u}
    s_{\CURL,h}(\uvec{u}_h,\uvec{u}_h)\to 0.
  \end{gather}

  Since $\Image\uG{k}{h}\subset\ker\uC{k}{h}$, the orthogonality assumption also gives $\uvec{u}_h\in(\Image\uG{k}{h})^\perp$.
  Using Remark \ref{rem:sobolev.potential} followed by %
  Theorem \ref{thm:sobolev.curl} with $s = 1$ (which is possible since $\Omega$ is convex) so that $s^* = 6$, %
  we have
  \begin{equation*}
    \norm{\bvec{L}^6(\Omega)}{\Pcurl{k}{h}\uvec{u}_h}
    \lesssim \tnorm{\CURL,h,6}{\uvec{u}_h}
    \lesssim \norm{\CURL,h}{\uvec{u}_h} + \norm{\DIV,h}{\uC{k}{h}\uvec{u}_h}
    \overset{\eqref{eq:conv.Ch.u.weak:basic-bounds}}\lesssim 1.
  \end{equation*}
  Since $(\Pcurl{k}{h}\uvec{u}_h)_{h\in\mathcal H}$ is uniformly bounded in $\bvec{L}^6(\Omega)$ and converges strongly to $\bvec{u}$ in $\bvec{L}^2(\Omega)$, weak compactness in $\bvec{L}^6(\Omega)$ and uniqueness of the limit yield $\bvec{u}\in\bvec{L}^6(\Omega)$.
  Interpolation between this uniform bound and \eqref{eq:conv.Pcurlh.u.L2} gives strong convergence in $\bvec{L}^r(\Omega)$ for every $r\in\lbrack 2,6\rparen$.
  The boundedness of $\Omega$ extends this convergence to every $r\in\lbrack 1,2\rparen$.
  This proves \eqref{eq:conv.Ph.u.strong}.
  \smallskip

  We next identify $\bvec{\xi}$. Let $\bvec{\phi}\in\vC{\infty}_c(\Omega)$.
  The primal consistency estimate of \cite[Theorem 8]{Chaumont-Frelet.Droniou.ea:25} yields
  \begin{multline*}
\left|
\int_\Omega \Pdiv{k}{h}\uC{k}{h}\uvec{u}_h\cdot\bvec{\phi}
-(\uC{k}{h}\uvec{u}_h,\qIdiv{k}{h}\bvec{\phi})_{\DIV,h}
\right|\\
  \lesssim \left(\norm{\bvec{L}^2(\Omega)}{\vlproj{k}{h}\bvec{\phi}-\bvec{\phi}}+h\norm{{L}^2(\Omega)}{\lproj{k}{h}(\DIV\bvec{\phi})-\DIV\bvec{\phi}}\right)\norm{\DIV,h}{\uC{k}{h}\uvec{u}_h}\to 0\text{ as $h\to 0$.}
  \end{multline*}
  On the other hand, the adjoint consistency estimate of \cite[Theorem 11(ii)]{Chaumont-Frelet.Droniou.ea:25} shows that
  \begin{align*}
  \Big|(\uC{k}{h}\uvec{u}_h,\qIdiv{k}{h}\bvec{\phi})_{\DIV,h}-{}&\int_\Omega\CURL\bvec{\phi}\cdot \Pcurl{k}{h}\uvec{u}_h\Big|\\
  \lesssim{}& \left(\norm{\bvec{L}^2(\Omega)}{\vlproj{k}{h}\bvec{\phi}-\bvec{\phi}}+h\norm{{L}^2(\Omega)}{\lproj{k}{h}(\DIV\bvec{\phi})-\DIV\bvec{\phi}}\right)\tnorm{\DIV,h,2}{\uC{k}{h}\uvec{u}_h}\\
  &+\norm{\bvec{L}^2(\Omega)}{\vlproj{k}{h}(\CURL\bvec{\phi})-\CURL\bvec{\phi}}\tnorm{\CURL,h,2}{\uvec{u}_h},
  \end{align*}
  which also tends to $0$ as $h\to 0$.
   Combining these two convergences gives
  \begin{equation*}
    \int_\Omega
    \Pdiv{k}{h}\uC{k}{h}\uvec{u}_h\cdot\bvec{\phi}
    -\int_\Omega \Pcurl{k}{h}\uvec{u}_h\cdot\CURL\bvec{\phi}
    \to0.
  \end{equation*}
  Indeed, the boundary term in the adjoint consistency relation vanishes since $\bvec{\phi}$ is compactly supported, while the consistency errors tend to zero since $\bvec{\phi}$ is smooth and the norms of $\uvec{u}_h$ and $\uC{k}{h}\uvec{u}_h$ are uniformly bounded.
  Passing to the limit using \eqref{eq:conv.Ch.u.xi} and \eqref{eq:conv.Pcurlh.u.L2}, we obtain
  \begin{equation*}
    \int_\Omega\bvec{\xi}\cdot\bvec{\phi}
    =\int_\Omega\bvec{u}\cdot\CURL\bvec{\phi}
    \qquad \forall\bvec{\phi}\in\vC{\infty}_c(\Omega).
  \end{equation*}
  Thus, $\CURL\bvec{u}=\bvec{\xi}$ in the sense of distributions. Since $\bvec{\xi}\in\bvec{L}^2(\Omega)$, this shows that $\bvec{u}\in\Hcurl{\Omega}$, and \eqref{eq:conv.Ch.u.xi} is precisely \eqref{eq:conv.Ch.u.weak}.
  \smallskip

  It remains to prove that $\bvec{u}\perp \ker\CURL$. Let $\bvec{w}\in\ker\CURL$. The cochain property of the quasi-interpolators \cite[Theorem 7]{Chaumont-Frelet.Droniou.ea:25} gives $\qIcurl{k}{h}\bvec{w}\in\ker\uC{k}{h}$.
  Therefore, the orthogonality of $\uvec{u}_h$ gives
  \begin{equation}\label{eq:orthogonality.uh.qIcurlw}
    \begin{aligned}
      0
      &= (\uvec{u}_h,\qIcurl{k}{h}\bvec{w})_{\CURL,h}\\
      \overset{\eqref{eq:def.L2.product.global}}&= \int_\Omega \Pcurl{k}{h}\uvec{u}_h\cdot\Pcurl{k}{h}\qIcurl{k}{h}\bvec{w}
      +s_{\CURL,h}(\uvec{u}_h,\qIcurl{k}{h}\bvec{w}).
    \end{aligned}
  \end{equation}
  The primal consistency properties and the convergence of the corresponding approximation errors \cite[Theorem 8 and Lemma 17]{Chaumont-Frelet.Droniou.ea:25} yield
  \begin{gather}\label{eq:conv.Pcurlh.qIcurlw}
    \Pcurl{k}{h}\qIcurl{k}{h}\bvec{w}
    \to\bvec{w} \quad \text{ strongly in $\bvec{L}^2(\Omega)$},\\
    \label{eq:conv.scurl.qIcurlw}
    s_{\CURL,h}(\qIcurl{k}{h}\bvec{w},\qIcurl{k}{h}\bvec{w}) \to 0.
  \end{gather}
  Hence, the Cauchy--Schwarz inequality for the positive semidefinite form $s_{\CURL,h}$, \eqref{eq:conv.stab.curlh.u} and \eqref{eq:conv.scurl.qIcurlw} give
  \begin{equation*}
    \begin{aligned}
      | s_{\CURL,h}(\uvec{u}_h,\qIcurl{k}{h}\bvec{w}) |
      &\le s_{\CURL,h}(\uvec{u}_h,\uvec{u}_h)^{\frac12} \,
      s_{\CURL,h} (\qIcurl{k}{h}\bvec{w},\qIcurl{k}{h}\bvec{w})^{\frac12} \to0.
    \end{aligned}
  \end{equation*}
  Passing to the limit in \eqref{eq:orthogonality.uh.qIcurlw} and using \eqref{eq:conv.Pcurlh.u.L2} and \eqref{eq:conv.Pcurlh.qIcurlw}, we conclude that $\int_\Omega\bvec{u}\cdot\bvec{w}=0$.
  Since $\bvec{w}\in\ker\CURL$ was arbitrary, we infer that $\bvec{u}\perp\ker\CURL$.
\end{proof}

\begin{lemma}[Approximation of functions in $\Hcurl{\Omega}\cap(\Image\GRAD)^\perp$]\label{lem:approximation.Hcurl}
  Assume that $\Omega$ is convex and let $(\Mh)_{h\in\mathcal H}$ be a sequence of quasi-uniform regular meshes.
  Let $\bvec{\phi}\in \Hcurl{\Omega}$ be such that $(\bvec{\phi},\GRAD q)_{L^2(\Omega)}=0$ for all $q\in \Hgrad{\Omega}$.
  Then, for all $h\in\mathcal H$, there exists $\wIcurl{k}{h}\bvec{\phi}\in\Xcurl{k}{\Th}$ such that
  \begin{alignat}{2}
    \label{eq:conv.tildeIphi}
    \norm{\bvec{L}^r(\Omega)}{\Pcurl{k}{h}(\wIcurl{k}{h}\bvec{\phi})-\bvec{\phi}}+
    s_{\CURL,h}(\wIcurl{k}{h}\bvec{\phi},\wIcurl{k}{h}\bvec{\phi})\to 0&\quad \text{as $h\to 0$, $\forall r\in [1,6)$},\\
    \label{eq:conv.curl.tildeIphi}
    \norm{\bvec{L}^2(\Omega)}{\Pdiv{k}{h}\uC{k}{h}(\wIcurl{k}{h}\bvec{\phi})-\CURL\bvec{\phi}}+
    s_{\DIV,h}(\uC{k}{h}(\wIcurl{k}{h}\bvec{\phi}),\uC{k}{h}(\wIcurl{k}{h}\bvec{\phi}))\to 0&\quad \text{as $h\to 0$},\\
    \label{eq:tildeIphi.orth}
    (\wIcurl{k}{h}\bvec{\phi},\uG{k}{h}\underline{q}_h)_{\CURL,h}=0&\quad\forall \underline{q}_h\in\Xgrad{k}{\Th}.
  \end{alignat}
\end{lemma}

\begin{proof}
  Define $\wIcurl{k}{h}\bvec{\phi}\in\Xcurl{k}{\Th}$ as the orthogonal projection of $\qIcurl{k}{h}\bvec{\phi}$ onto
  $(\Image\uG{k}{h})^\perp$. The relation \eqref{eq:tildeIphi.orth} holds by construction and, moreover, there exists $\underline{q}_h\in\Xgrad{k}{\Th}$ such that
  \begin{equation}\label{eq:diff.interp}
    \qIcurl{k}{h}\bvec{\phi}-\wIcurl{k}{h}\bvec{\phi}=\uG{k}{h}\underline{q}_h.
  \end{equation}
  The complex property $\uC{k}{h}\circ\uG{k}{h}=\uvec{0}_h$ shows that $\uC{k}{h}\wIcurl{k}{h}\bvec{\phi}=\uC{k}{h}\qIcurl{k}{h}\bvec{\phi}=\qIdiv{k}{h}(\CURL\bvec{\phi})$, where the last equality follows from the commutation properties of the quasi-interpolators. By the primal consistency of $\Pdiv{k}{h}$ and $\qIdiv{k}{h}$ stated in \cite[Theorem 8]{Chaumont-Frelet.Droniou.ea:25}, we deduce that \eqref{eq:conv.curl.tildeIphi} holds.
  \smallskip

  To establish \eqref{eq:conv.tildeIphi}, we start by writing
  \begin{align*}
    \norm{\CURL,h}{\qIcurl{k}{h}\bvec{\phi}-\wIcurl{k}{h}\bvec{\phi}}^2\overset{\eqref{eq:diff.interp}}&=
    (\qIcurl{k}{h}\bvec{\phi}-\wIcurl{k}{h}\bvec{\phi},\uG{k}{h}\underline{q}_h)_{\CURL,h}\\
    \overset{\eqref{eq:tildeIphi.orth}}&=(\qIcurl{k}{h}\bvec{\phi},\uG{k}{h}\underline{q}_h)_{\CURL,h}\\
    &\lesssim \left(
    \norm{\bvec{L}^2(\Omega)}{\vlproj{k}{h}\bvec{\phi}-\bvec{\phi}}
    + h \norm{\bvec{L}^2(\Omega)}{\vlproj{k}{h}(\CURL\bvec{\phi})-\CURL\bvec{\phi}} 
    \right)
    \norm{\CURL,h}{\uG{k}{h}\underline{q}_h}
  \end{align*}
  where the conclusion follows, as in \eqref{eq:use.adjoint}, from the adjoint consistency stated in \cite[Theorem 11, (i)]{Chaumont-Frelet.Droniou.ea:25} together with $\DIV\bvec{\phi}=0$ and $\bvec{\phi}\cdot\normal_\Omega=0$ on $\partial\Omega$ (see \eqref{eq:im.grad.div.0}). Recalling \eqref{eq:diff.interp}, simplifying and using the definition of $\norm{\CURL,h}{{\cdot}}$ stemming from \eqref{eq:def.L2.product.global} shows that
  \begin{equation*}
    \norm{\bvec{L}^2(\Omega)}{\Pcurl{k}{h}(\qIcurl{k}{h}\bvec{\phi})-\Pcurl{k}{h}(\wIcurl{k}{h}\bvec{\phi})}
    +s_{\CURL,h}(\qIcurl{k}{h}\bvec{\phi}-\wIcurl{k}{h}\bvec{\phi},\qIcurl{k}{h}\bvec{\phi}-\wIcurl{k}{h}\bvec{\phi})\to 0\text{ as $h\to 0$.}
  \end{equation*}
  Using \eqref{eq:conv.Pcurlh.qIcurlw} and \eqref{eq:conv.scurl.qIcurlw} (with $\bvec{w}=\bvec{\phi}$) together with triangle inequalities, we infer that \eqref{eq:conv.tildeIphi} holds with $r=2$.

  By \eqref{eq:qIbullet.bound}, $\qIcurl{k}{h}\bvec{\phi}$ and $\qIdiv{k}{h}(\CURL\bvec{\phi})$ are bounded, respectively, in $\Xcurl{k}{\Th}$ and $\Xdiv{k}{\Th}$. Since $\wIcurl{k}{h}\bvec{\phi}$ is an orthogonal projection of $\qIcurl{k}{h}\bvec{\phi}$ and
  $\uC{k}{h}\wIcurl{k}{h}\bvec{\phi}=\qIdiv{k}{h}(\CURL\bvec{\phi})$, we infer that $\norm{\CURL,h}{\wIcurl{k}{h}\bvec{\phi}}+\norm{\DIV,h}{\uC{k}{h}\wIcurl{k}{h}\bvec{\phi}}$ is bounded. Recalling \eqref{eq:tildeIphi.orth}, the discrete Sobolev embedding of Theorem \ref{thm:sobolev.curl} therefore shows that $\Pcurl{k}{h}(\wIcurl{k}{h}\bvec{\phi})$ is bounded in $\bvec{L}^6(\Omega)$.
  Combined with \eqref{eq:conv.tildeIphi} for $r=2$, this proves that this convergence holds for all $r<6$. \end{proof}

\begin{lemma}[Inverse estimates]
  If $\Mh$ is a regular mesh, then, for all $T\in\Th$,
  \begin{equation}
    \label{eq:bound.L4.PqI}
    \norm{\bvec{L}^4(T)}{\Pcurl{k}{T}\qIcurl{k}{T}\bvec{w}}
    \lesssim h_T^{-\frac34}\left(\norm{\bvec{L}^2(\widetilde{T})}{\bvec{w}}+h_T\norm{\bvec{L}^2(\widetilde{T})}{\CURL\bvec{w}}\right)\quad\forall \bvec{w}\in \Hcurl{\widetilde{T}},
  \end{equation}
  where we recall that $\widetilde{T}$ is the patch of cells around $T$.
  If $\Mh$ is moreover quasi-uniform, then
  \begin{equation}
    \label{eq:inverse.est.rh}
    \norm{\CURL,h}{\uG{k}{h}\underline{r}_h}\lesssim h^{-\frac34}\norm{P}{\underline{r}_h}\qquad\forall \underline{r}_h\in\Xgrad{k}{\Th}.
  \end{equation}
\end{lemma}

\begin{proof}
  Using the discrete Lebesgue inequality for polynomials \cite[Lemma 1.25]{Di-Pietro.Droniou:20}, we have
  \[
  \norm{\bvec{L}^4(T)}{\Pcurl{k}{T}\qIcurl{k}{T}\bvec{w}}\lesssim h_T^{-\frac34}\norm{\bvec{L}^2(T)}{\Pcurl{k}{T}\qIcurl{k}{T}\bvec{w}}\lesssim h_T^{-\frac34}\tnorm{\CURL,T,2}{\qIcurl{k}{T}\bvec{w}},
  \]
  where the last inequality follows from \eqref{eq:local.norm.equivalence} and the definition of $\norm{\CURL,T}{{\cdot}}$. Invoking the bound \eqref{eq:qIbullet.bound} on the local quasi-interpolator concludes the proof of \eqref{eq:bound.L4.PqI}.
  \smallskip

  Now let $\underline{r}_h\in\Xgrad{k}{\Th}$. We have
  \begin{equation} \label{eq:bound.uGr.hneg}
    \norm{\CURL, h}{\uG{k}{h}\underline{r}_h}^2=(\uG{k}{h}\underline{r}_h,\uG{k}{h}\underline{r}_h)_{\CURL,h}
    \overset{\eqref{eq:def.norm.P}}\lesssim\norm{P}{\underline{r}_h}\tnorm{\CURL,h,4}{\uG{k}{h}\underline{r}_h}.
  \end{equation}
  The discrete Lebesgue inequality \eqref{eq:discrete.lebesgue}, applied to each $T\in\Th$ with $(\alpha,\beta)=(4,2)$, and the quasi-uniformity of the mesh yield
  \[
  \tnorm{\CURL,h,4}{\uG{k}{h}\underline{r}_h}
  = \left(\sum_{T\in\Th}\tnorm{\CURL,T,4}{\uG{k}{T}\underline{r}_T}^4\right)^{\frac14}
  \!\lesssim h^{-\frac34}\left(\sum_{T\in\Th}\big(\tnorm{\CURL,T,2}{\uG{k}{T}\underline{r}_T}^2\big)^2\right)^{\frac14}
  \overset{\eqref{eq:sum.xi.3}}\lesssim h^{-\frac34}\tnorm{\CURL,h,2}{\uG{k}{h}\underline{r}_h}.
  \]
  Plugging this estimate into \eqref{eq:bound.uGr.hneg}, invoking the norm equivalence \eqref{eq:global.norm.equivalence} and simplifying concludes the proof of \eqref{eq:inverse.est.rh}.
\end{proof}

\begin{lemma}[Consistency of $\Icurl{k}{h}$ and $\qIcurl{k}{h}$ in non-Hilbertian settings]
  If $(\Mh)_{h\in\mathcal H}$ is a sequence of regular meshes, then
  \begin{alignat}{2}
    \label{eq:conv.PI.Linfty}
    \norm{\bvec{L}^\infty(\Omega)}{\Pcurl{k}{h}(\Icurl{k}{h}\bvec{v})-\bvec{v}}
    +s_{\CURL,h}(\Icurl{k}{h}\bvec{v},\Icurl{k}{h}\bvec{v})&\to 0\text{ as $h\to 0$}, &&\quad\forall \bvec{v}\in \bvec{C}(\overline{\Omega}),\\
    \label{eq:conv.PqI.L4}
    \norm{\bvec{L}^4(\Omega)}{\Pcurl{k}{h}(\qIcurl{k}{h}\bvec{w})-\bvec{w}}&\lesssim h^{\frac14}\norm{\bvec{H}^1(\Omega)}{\bvec{w}}
    &&\quad \forall \bvec{w}\in \bvec{H}^1(\Omega).
  \end{alignat}
\end{lemma}

\begin{proof}
  Let $\bvec{v}\in \bvec{C}(\overline{\Omega})$.
  By polynomial consistency of $\Pcurl{k}{h}\circ\Icurl{k}{h}$ (see \cite[Eq.~(4.7)]{Di-Pietro.Droniou:23}) we have, for all $T\in\Th$,
  \[
  \Pcurl{k}{T}(\Icurl{k}{T}\bvec{v})-\bvec{v}
  = \Pcurl{k}{T}(\Icurl{k}{T}(\bvec{v}-\vlproj{k}{T}\bvec{v}))-(\bvec{v}-\vlproj{k}{T}\bvec{v})
  \]
  and thus
  \begin{align}
    \norm{\bvec{L}^\infty(\Omega)}{\Pcurl{k}{h}(\Icurl{k}{h}\bvec{v})-\bvec{v}}
    &=\max_{T\in\Th}\norm{\bvec{L}^\infty(T)}{\Pcurl{k}{T}(\Icurl{k}{T}(\bvec{v}-\vlproj{k}{T}\bvec{v}))-(\bvec{v}-\vlproj{k}{T}\bvec{v})}\nonumber\\
    &\lesssim \max_{T\in\Th}\norm{\bvec{L}^\infty(T)}{\bvec{v}-\vlproj{k}{T}\bvec{v}}
    \label{eq:est.PIv-v}
  \end{align}
  where we have used the bound on $\Pcurl{k}{h}$ established in \cite[Eq.~(A.10)]{Di-Pietro.Droniou.ea:24} (with $s=\infty$) and the fact that the $L^\infty$ norms of the components of $\Icurl{k}{T}(\bvec{v}-\vlproj{k}{T}\bvec{v})$ are bounded by $\norm{\bvec{L}^\infty(T)}{\bvec{v}-\vlproj{k}{T}\bvec{v}}$ (which follows from $L^\infty$-boundedness of $L^2$-projectors on polynomial spaces \cite[Lemma 1.44]{Di-Pietro.Droniou:20}). Using first the polynomial consistency of $s_{\CURL,T}$ and then the local norm equivalence \eqref{eq:local.norm.equivalence}, we have
  \begin{align}
    s_{\CURL,T}(\Icurl{k}{T}\bvec{v},\Icurl{k}{T}\bvec{v})
    ={}&s_{\CURL,T}(\Icurl{k}{T}(\bvec{v}-\vlproj{k}{T}\bvec{v}),\Icurl{k}{T}(\bvec{v}-\vlproj{k}{T}\bvec{v}))\nonumber\\
    \lesssim{}&
    \tnorm{\CURL,T,2}{\Icurl{k}{h}(\bvec{v}-\vlproj{k}{T}\bvec{v})}^{2}
    \nonumber\\
    \lesssim{}& |T|\norm{\bvec{L}^\infty(T)}{\bvec{v}-\vlproj{k}{T}\bvec{v}}^{2},
    \label{eq:bound.sT.max}
  \end{align}
  where the last inequality follows from the bound given above on the components of $\Icurl{k}{h}(\bvec{v}-\vlproj{k}{T}\bvec{v})$.
  Summing \eqref{eq:bound.sT.max} over $T\in\Th$, applying H\"older's inequality to the sum with exponents $(1,\infty)$ and adding \eqref{eq:est.PIv-v} yields
  \[
  \norm{\bvec{L}^\infty(\Omega)}{\Pcurl{k}{h}(\Icurl{k}{h}\bvec{v})-\bvec{v}}
  + s_{\CURL,h}(\Icurl{k}{h}\bvec{v},\Icurl{k}{h}\bvec{v})^{\frac12}    
  \lesssim \max_{T\in\Th}\norm{\bvec{L}^\infty(T)}{\bvec{v}-\vlproj{k}{T}\bvec{v}}.
  \]
  Since $\bvec{v}$ is uniformly continuous on $\overline{\Omega}$, the right-hand side above tends to $0$ as $h\to 0$, which establishes \eqref{eq:conv.PI.Linfty}.
  \smallskip

  Now take $\bvec{w}\in\bvec{H}^1(\Omega)$. For all $T\in\Th$, the potential reconstruction and quasi-interpolator satisfy the polynomial consistency property $\Pcurl{k}{T}(\qIcurl{k}{T}\bvec{z})=\bvec{z}$ for all $\bvec{z}\in\vPoly{k}(\widetilde{T})$, the space of vector-valued polynomials of total degree $\le k$ on the patch $\widetilde{T}$.
  Hence, with $\vlproj{k}{\widetilde{T}}$ denoting the $L^2$-orthogonal projector onto $\vPoly{k}(\widetilde{T})$, the triangle inequality gives
  \begin{equation}
    \begin{aligned}
      &\|\Pcurl{k}{T}(\qIcurl{k}{T}\bvec{w})-\bvec{w}\|_{\bvec{L}^4(T)}
      \\
      &\quad
      \begin{aligned}[t]
        &= \norm{\bvec{L}^4(T)}{\Pcurl{k}{T}(\qIcurl{k}{T}(\bvec{w}-\vlproj{k}{\widetilde{T}}\bvec{w}))-(\bvec{w}-\vlproj{k}{\widetilde{T}}\bvec{w})}
        \\
        \overset{\eqref{eq:bound.L4.PqI}}&\lesssim
        h_T^{-\frac34}\left(\norm{\bvec{L}^2(\widetilde{T})}{\bvec{w}-\vlproj{k}{\widetilde{T}}\bvec{w}}+h_T\norm{\bvec{L}^2(\widetilde{T})}{\CURL(\bvec{w}-\vlproj{k}{\widetilde{T}}\bvec{w})}\right)
        +\norm{\bvec{L}^4(T)}{\bvec{w}-\vlproj{k}{\widetilde{T}}\bvec{w}}.
      \end{aligned}
    \end{aligned}
    \label{eq:approx.PqI.optimal}
  \end{equation}
  An interpolation estimate gives
  \[
\norm{\bvec{L}^4(T)}{\bvec{w}-\vlproj{k}{\widetilde{T}}\bvec{w}}
    \le \norm{\bvec{L}^2(T)}{\bvec{w}-\vlproj{k}{\widetilde{T}}\bvec{w}}^{\frac14}\norm{\bvec{L}^6(T)}{\bvec{w}-\vlproj{k}{\widetilde{T}}\bvec{w}}^{\frac34}
    \lesssim \norm{\bvec{L}^2(T)}{\bvec{w}-\vlproj{k}{\widetilde{T}}\bvec{w}}^{\frac14}\norm{\bvec{L}^6(\widetilde{T})}{\bvec{w}}^{\frac34},
  \]
  where the conclusion follows from $L^6$-boundedness of $\vlproj{k}{\widetilde{T}}$ (see \cite[Lemma 1.44]{Di-Pietro.Droniou:20} and note that $\widetilde{T}$ is connected by star-shaped sets in the sense of \cite[Definition 1.41]{Di-Pietro.Droniou:20}).
  Plugging this estimate into \eqref{eq:approx.PqI.optimal} and using the approximation properties and $H^1$-boundedness of $\vlproj{k}{\widetilde{T}}$, we infer
  \[
  \norm{\bvec{L}^4(T)}{\Pcurl{k}{T}(\qIcurl{k}{T}\bvec{w})-\bvec{w}}
  \lesssim h_T^{\frac14}\seminorm{\bvec{H}^1(\widetilde{T})}{\bvec{w}}+h_T^{\frac14}\seminorm{\bvec{H}^1(\widetilde{T})}{\bvec{w}}^{\frac14}\norm{\bvec{L}^6(\widetilde{T})}{\bvec{w}}^{\frac34}.
  \]
  Taking the discrete $\ell^4(\Th)$ norm of this inequality and applying the triangle inequality gives
  \begin{align*}
    \norm{\bvec{L}^4(\Omega)}{\Pcurl{k}{h}(\qIcurl{k}{h}\bvec{w})-\bvec{w}}
    &\lesssim h^{\frac14}\left[
      \sum_{T\in\Th}\left(
      \seminorm{\bvec{H}^1(\widetilde{T})}{\bvec{w}}^2
      \right)^2
    \right]^{\frac14}
    + h^{\frac14}\left(\sum_{T\in\Th}\seminorm{\bvec{H}^1(\widetilde{T})}{\bvec{w}}\norm{\bvec{L}^6(\widetilde{T})}{\bvec{w}}^{3}\right)^{\frac14}\\
    \overset{\eqref{eq:sum.xi.3},\eqref{eq:norm.patches.Lrho}}&\lesssim h^{\frac14}\seminorm{\bvec{H}^1(\Omega)}{\bvec{w}}
    +h^{\frac14}\seminorm{\bvec{H}^1(\Omega)}{\bvec{w}}^{\frac14}\norm{\bvec{L}^6(\Omega)}{\bvec{w}}^{\frac34},
  \end{align*}
  where we have additionally applied the Cauchy--Schwarz inequality to the last sum in the first inequality.
  Applying the Sobolev embedding $\bvec{H}^1(\Omega)\hookrightarrow \bvec{L}^6(\Omega)$ completes the proof of \eqref{eq:conv.PqI.L4}.
\end{proof}

\begin{remark}[Optimal approximation properties of $\Pcurl{k}{h}\circ\qIcurl{k}{h}$]
  The bound \eqref{eq:conv.PqI.L4} is sufficient to prove the convergence of the pressure in Theorem \ref{thm:convergence.compactness}. The arguments above can, however, be used to obtain an optimal approximation property of $\Pcurl{k}{h}\circ\qIcurl{k}{h}$ in the $L^4$-norm for more regular functions. Indeed, applying H\"older's inequality to the $L^2$-norms in \eqref{eq:approx.PqI.optimal}, we obtain
  \[
  \begin{aligned}
    \norm{\bvec{L}^4(T)}{\Pcurl{k}{T}(\qIcurl{k}{T}\bvec{w})-\bvec{w}}
    &\lesssim
    \norm{\bvec{L}^4(\widetilde{T})}{\bvec{w}-\vlproj{k}{\widetilde{T}}\bvec{w}}
    +h_T\norm{\bvec{L}^4(\widetilde{T})}{\CURL(\bvec{w}-\vlproj{k}{\widetilde{T}}\bvec{w})}
    \\
    &\lesssim h_T^{k+1}\seminorm{\bvec{W}^{k+1,4}(\widetilde{T})}{\bvec{w}},
  \end{aligned}
  \]
  where the last inequality follows from the approximation properties of $\vlproj{k}{\widetilde{T}}$.
  The same arguments can be applied to any $L^s$-norm with $s\ge 2$ (not just $s=4$).

  Note that optimal (and sharper) $L^2$-approximation properties have been obtained for $\Pcurl{k}{h}\circ\qIcurl{k}{h}$ in \cite{Chaumont-Frelet.Droniou.ea:25}, and that the results in the appendix of \cite{Di-Pietro.Droniou.ea:24} give optimal $L^s$-approximation properties for $\Pcurl{k}{h}\circ\Icurl{k}{h}$ for any $s\in [1,\infty]$.
\end{remark}

\subsubsection{Proof of Theorem \ref{thm:convergence.compactness}}\label{sec:proof.convergence.compactness}

The proof of the convergence of the DDR scheme \eqref{eq:ddr.scheme} to a solution of the continuous model \eqref{eq:ns.curl} proceeds in several steps.
\medskip\\
\emph{Step 1: bounds.}
The definition of the DDR interpolator gives
\begin{equation}\label{eq:bound.If.compactness}
  \norm{\CURL,h}{\Icurl{k}{h}\bvec{f}}
  \lesssim \norm{\bvec{C}(\overline{\Omega})}{\bvec{f}}
  \lesssim 1.
\end{equation}
Since $\Omega$ is convex, its topology is trivial and the exactness and cohomology properties of the DDR complex \cite{Di-Pietro.Droniou:23,Di-Pietro.Droniou.ea:23} yield $\ker\uC{k}{h}=\Image\uG{k}{h}$.
The relation \eqref{eq:ddr.scheme.incompressibility} therefore implies $\uvec{u}_h\perp\ker\uC{k}{h}$.
Again using the convexity of $\Omega$, we see that Theorem \ref{thm:sobolev.curl} applies with $s=1$ and $s^*=6$.
Combining this theorem with Remark \ref{rem:sobolev.potential}, the discrete Poincar\'e inequality \cite{Di-Pietro.Droniou.ea:25*1}, and the continuous embedding $L^6(\Omega)\hookrightarrow L^4(\Omega)$ gives
\begin{equation*}
  \norm{\bvec{L}^4(\Omega)}{\Pcurl{k}{h}\uvec{v}_h}
  \lesssim \norm{\DIV,h}{\uC{k}{h}\uvec{v}_h}
  \qquad\forall\uvec{v}_h\in(\ker\uC{k}{h})^\perp.
\end{equation*}
Consequently, the discrete Sobolev constants appearing in \cite[Lemma~5]{Di-Pietro.Droniou.ea:24} are uniformly bounded.
Recalling the a priori estimates of \cite[Lemmas~4 and~5]{Di-Pietro.Droniou.ea:24} and using \eqref{eq:bound.If.compactness}, we obtain
\begin{equation}\label{eq:apriori.compactness}
  \norm{\DIV,h}{\uC{k}{h}\uvec{u}_h}
  +\norm{\CURL,h}{\uvec{u}_h}
  +\norm{P}{\underline{p}_h}
  \lesssim 1.
\end{equation}

\medskip
\noindent\emph{Step 2: convergence of the velocity.}
By \eqref{eq:apriori.compactness} and Lemma \ref{lem:compactness.curl}, there exists $\bvec{u}\in\Hcurl{\Omega}$ such that $\bvec{u}\perp\ker\CURL$ and, up to a subsequence (not relabelled), \eqref{eq:conv.Ph.u.strong} and \eqref{eq:conv.Ch.u.weak} hold, which establishes \eqref{eq:conv.Pcurlh.u}.
The relation $\bvec{u}\perp\ker\CURL$ gives \eqref{eq:ns.curl.incompressibility}, and it remains to prove that the momentum equation holds.

In this step, we only prove that \eqref{eq:ns.curl.momentum} holds for test functions $\bvec{\phi}\in\Hcurl{\Omega}\cap (\Image\GRAD)^\perp$, which eliminates the pressure term; note that, by \eqref{eq:im.grad.div.0} and \eqref{eq:regularity.Hcurl.Hdiv}, such functions $\bvec{\phi}$ belong to $\bvec{L}^4(\Omega)$. Recalling the interpolator $\wIcurl{k}{h}$ given by Lemma \ref{lem:approximation.Hcurl}, we let $\uvec{v}_h=\wIcurl{k}{h}\bvec{\phi}$ in \eqref{eq:ddr.scheme.momentum} and use \eqref{eq:tildeIphi.orth} together with the definition of $(\cdot,\cdot)_{\CURL,h}$ to get
\begin{equation} \label{eq:plug.in.scheme}
  \begin{aligned}
    \nu\int_\Omega\Pdiv{k}{h}\uC{k}{h}\uvec{u}_h\cdot&\Pdiv{k}{h}\uC{k}{h}(\wIcurl{k}{h}\bvec{\phi})
    +\nu s_{\DIV,h}(\uC{k}{h}\uvec{u}_h,\uC{k}{h}\wIcurl{k}{h}\bvec{\phi})\\
    &+\int_\Omega [(\Pdiv{k}{h}\uC{k}{h}\uvec{u}_h)\times\Pcurl{k}{h}\uvec{u}_h]\cdot\Pcurl{k}{h} (\wIcurl{k}{h}\bvec{\phi})
    =(\Icurl{k}{h}\bvec{f},\wIcurl{k}{h}\bvec{\phi})_{\CURL,h}.
  \end{aligned}
\end{equation}
The weak and strong convergence properties \eqref{eq:conv.Ch.u.weak} and \eqref{eq:conv.curl.tildeIphi}, the Cauchy--Schwarz inequality for $s_{\DIV,h}$ and the bound \eqref{eq:apriori.compactness} allow us to pass to the limit in the first two terms:
\begin{equation}\label{eq:conv.velocity.1}
  \int_\Omega\Pdiv{k}{h}\uC{k}{h}\uvec{u}_h\cdot\Pdiv{k}{h}\uC{k}{h}(\wIcurl{k}{h}\bvec{\phi})
  +s_{\DIV,h}(\uC{k}{h}\uvec{u}_h,\uC{k}{h}\wIcurl{k}{h}\bvec{\phi})\to
  \int_\Omega \CURL\bvec{u}\cdot \CURL\bvec{\phi}\quad\text{ as $h\to 0$}.
\end{equation}
Similarly, expanding the definition of $(\cdot,\cdot)_{\CURL,h}$ and using the convergence properties \eqref{eq:conv.tildeIphi} and \eqref{eq:conv.PI.Linfty}, we obtain
\begin{equation}\label{eq:conv.velocity.2}
  (\Icurl{k}{h}\bvec{f},\wIcurl{k}{h}\bvec{\phi})_{\CURL,h}\to
  \int_\Omega \bvec{f}\cdot \bvec{\phi}.
\end{equation}
To pass to the limit in the nonlinear term, we use \eqref{eq:conv.Ch.u.weak}, \eqref{eq:conv.Pcurlh.u} and \eqref{eq:conv.tildeIphi} with $r=4$:
\begin{equation}\label{eq:conv.velocity.3}
  \int_\Omega [(\Pdiv{k}{h}\uC{k}{h}\uvec{u}_h)\times\Pcurl{k}{h}\uvec{u}_h]\cdot\Pcurl{k}{h} (\wIcurl{k}{h}\bvec{\phi})
  \to \int_\Omega [(\CURL\bvec{u})\times \bvec{u}]\cdot \bvec{\phi}.
\end{equation}
Plugging \eqref{eq:conv.velocity.1}, \eqref{eq:conv.velocity.2} and \eqref{eq:conv.velocity.3} into \eqref{eq:plug.in.scheme} shows that
\begin{equation}\label{eq:lim.problem.u}
  \nu\int_\Omega \CURL\bvec{u}\cdot \CURL\bvec{\phi}+\int_\Omega [(\CURL\bvec{u})\times \bvec{u}]\cdot \bvec{\phi}=\int_\Omega \bvec{f}\cdot \bvec{\phi}\qquad\forall\bvec{\phi}\in \Hcurl{\Omega}\cap (\Image\GRAD)^\perp.
\end{equation}

To conclude this step, we prove \eqref{eq:conv.Ch.u} and \eqref{eq:conv.stab.Ch.u}, the strong convergence properties of the discrete curl. Take $\uvec{v}_h=\uvec{u}_h$ in \eqref{eq:ddr.scheme.momentum}. Eliminating the pressure contribution using \eqref{eq:ddr.scheme.incompressibility} and the nonlinear term using $[(\Pdiv{k}{h}\uC{k}{h}\uvec{u}_h)\times\Pcurl{k}{h}\uvec{u}_h]\cdot\Pcurl{k}{h}\uvec{u}_h=0$, we find
\[
\nu\norm{\bvec{L}^2(\Omega)}{\Pdiv{k}{h}\uC{k}{h}\uvec{u}_h}^2
+\nu s_{\DIV,h}(\uC{k}{h}\uvec{u}_h,\uC{k}{h}\uvec{u}_h)
= \nu(\uC{k}{h}\uvec{u}_h,\uC{k}{h}\uvec{u}_h)_{\DIV,h}
=(\Icurl{k}{h}\bvec{f},\uvec{u}_h)_{\CURL,h}.
\]
By the convergence properties \eqref{eq:conv.Ph.u.strong} and \eqref{eq:conv.PI.Linfty} and the bound on $s_{\CURL,h}(\uvec{u}_h,\uvec{u}_h)$ resulting from \eqref{eq:apriori.compactness}, we can pass to the limit in the right-hand side to get
\[
\nu \limsup_{h\to 0}\left(\norm{\bvec{L}^2(\Omega)}{\Pdiv{k}{h}\uC{k}{h}\uvec{u}_h}^2
+s_{\DIV,h}(\uC{k}{h}\uvec{u}_h,\uC{k}{h}\uvec{u}_h)\right)\le \int_\Omega \bvec{f}\cdot\bvec{u}.
\]
Since $\bvec{u}\in\Hcurl{\Omega}\cap (\Image\GRAD)^\perp$, we can use it as a test function in \eqref{eq:lim.problem.u} to substitute $\int_\Omega \bvec{f}\cdot\bvec{u}$ and obtain
\begin{equation}\label{eq:conv.Ch.strong.2}
  \nu \limsup_{h\to 0}\left(\norm{\bvec{L}^2(\Omega)}{\Pdiv{k}{h}\uC{k}{h}\uvec{u}_h}^2
  +s_{\DIV,h}(\uC{k}{h}\uvec{u}_h,\uC{k}{h}\uvec{u}_h)\right)\le \nu\norm{\bvec{L}^2(\Omega)}{\CURL\bvec{u}}^2.
\end{equation}
This shows that $\limsup_{h\to 0}\norm{\bvec{L}^2(\Omega)}{\Pdiv{k}{h}\uC{k}{h}\uvec{u}_h}\le \norm{\bvec{L}^2(\Omega)}{\CURL\bvec{u}}$, which, combined with the weak convergence \eqref{eq:conv.Ch.u.weak}, proves the strong convergence \eqref{eq:conv.Ch.u}. Using this strong convergence in \eqref{eq:conv.Ch.strong.2}, we obtain $\limsup_{h\to 0}s_{\DIV,h}(\uC{k}{h}\uvec{u}_h,\uC{k}{h}\uvec{u}_h)\le 0$, which establishes \eqref{eq:conv.stab.Ch.u}.

\medskip
\noindent\emph{Step 3: existence of the pressure.}

Since $\bvec{u}\in\Hcurl{\Omega}\cap \bvec{L}^4(\Omega)$, we have $\bvec{f}-(\CURL\bvec{u})\times \bvec{u}\in\bvec{L}^{\frac{4}{3}}(\Omega)$. By convexity of $\Omega$, \cite[Lemma 2.1]{Girault:90} gives the existence of a zero-average $p\in W^{1,\frac{4}{3}}(\Omega)$ such that
\begin{equation}\label{eq:eq.gradp}
  \int_\Omega \GRAD p\cdot\GRAD q=\int_\Omega (\bvec{f}-(\CURL\bvec{u})\times \bvec{u})\cdot\GRAD q
  \qquad\forall q\in W^{1,4}(\Omega).
\end{equation}
We can now check that $(\bvec{u},p)\in\Hcurl{\Omega}\times  W^{1,\frac{4}{3}}(\Omega)$ satisfies \eqref{eq:ns.curl}. The relation \eqref{eq:ns.curl.incompressibility} on $\bvec{u}$ has already been established. Let $\bvec{v}\in\Hcurl{\Omega}\cap \bvec{L}^4(\Omega)$ and consider its Hodge decomposition
\[
\bvec{v}=\bvec{\phi}+\GRAD q\quad\text{ with }\quad\bvec{\phi}\in \Hcurl{\Omega}\cap(\Image\GRAD)^\perp\,,
\quad q\in W^{1,4}(\Omega).
\]
The regularity of $q$ beyond $\Hgrad{\Omega}$ follows from the observation that, by convexity, $\bvec{\phi}\in\bvec{H}^1(\Omega)\subset\bvec{L}^4(\Omega)$, so $\GRAD q=\bvec{v}-\bvec{\phi}\in\bvec{L}^4(\Omega)$.
Apply \eqref{eq:lim.problem.u} to $\bvec{\phi}$ and \eqref{eq:eq.gradp} to $q$, then sum the two equations to obtain
\begin{align*}
  \nu\int_\Omega \CURL\bvec{u}\cdot \CURL\bvec{\phi}+\int_\Omega [(\CURL\bvec{u})\times \bvec{u}]\cdot (\bvec{\phi}+\GRAD q)+\int_\Omega \GRAD p\cdot \GRAD q=\int_\Omega \bvec{f}\cdot (\bvec{\phi}+\GRAD q).
\end{align*}
Since $\CURL\circ\GRAD=0$, we have $ \CURL\bvec{\phi}=\CURL(\bvec{\phi}+\GRAD q)=\CURL\bvec{v}$ and, since $\bvec{\phi}\perp\Image\GRAD$, we note that
$\int_\Omega \GRAD p\cdot \GRAD q=\int_\Omega \GRAD p\cdot (\bvec{\phi}+\GRAD q)=\int_\Omega\GRAD p\cdot\bvec{v}$ 
(the orthogonality of $\bvec{\phi}$ to gradients extends by density from potentials in $H^1(\Omega)$ to potentials in $W^{1,\frac{4}{3}}(\Omega)$, since $\bvec{\phi}\in\bvec{L}^4(\Omega)$).
The above equation therefore proves that \eqref{eq:ns.curl.momentum} holds for $\bvec{v}$, which concludes this step.

\medskip
\noindent\emph{Step 4: convergence of the pressure.}

It remains to prove \eqref{eq:conv.Pgradh.p}, which relies on a classical duality argument.
Let $g\in L^2(\Omega)$ be such that $\int_\Omega g=0$ and let $\varphi\in H^1(\Omega)$ be the zero-average solution of $\Delta\varphi=g$ in $\Omega$, $\GRAD\varphi\cdot\normal_{\Omega}=0$ on $\partial\Omega$. Since $\Omega$ is convex, we have $\varphi\in H^2(\Omega)$ with
\begin{equation}\label{eq:reg.varphi}
  \norm{H^2(\Omega)}{\varphi}\lesssim\norm{L^2(\Omega)}{g}.
\end{equation}
We then write, introducing $\pm(\uG{k}{h}\underline{p}_h,\qIcurl{k}{h}(\GRAD\varphi))_{\CURL,h}$,
\begin{align}
  \int_\Omega (\Pgrad{k+1}{h}\underline{p}_h-p)g=&\left(\int_\Omega \Pgrad{k+1}{h}\underline{p}_h\DIV(\GRAD\varphi)
  +(\uG{k}{h}\underline{p}_h,\qIcurl{k}{h}(\GRAD\varphi))_{\CURL,h}\right)\nonumber\\
  &+\left(-(\uG{k}{h}\underline{p}_h,\qIcurl{k}{h}(\GRAD\varphi))_{\CURL,h}+\int_\Omega \GRAD p\cdot\GRAD\varphi\right)\nonumber\\
  ={}&\term_1+\term_2.
  \label{eq:conv.pressure.1}
\end{align}
The first term can be estimated using the adjoint consistency result in \cite[Theorem 11, (i)]{Chaumont-Frelet.Droniou.ea:25}: since $\GRAD\varphi\cdot\normal_\Omega=0$, we obtain
\begin{align}
  |\term_1|&\lesssim \left(\norm{\bvec{L}^2(\Omega)}{\GRAD\varphi-\vlproj{k}{h}(\GRAD\varphi)} + h\norm{L^2(\Omega)}{g-\lproj{k}{h}g}\right)\tnorm{\CURL,h,2}{\uG{k}{h}\underline{p}_h}\nonumber\\
  \overset{\eqref{eq:global.norm.equivalence}, \eqref{eq:inverse.est.rh}}&\lesssim h^{\frac14}\left(\norm{\bvec{H}^1(\Omega)}{\GRAD\varphi} + \norm{L^2(\Omega)}{g}\right)\norm{P}{\underline{p}_h}\nonumber\\
  \overset{\eqref{eq:reg.varphi},\eqref{eq:apriori.compactness}}&\lesssim h^{\frac14}\norm{L^2(\Omega)}{g},
  \label{eq:conv.pressure.2}
\end{align}
where, in the second inequality, we have moreover used the approximation properties of $\vlproj{k}{h}$ on $\GRAD\varphi$.

To estimate $\term_2$, we use the discrete and variational momentum equations \eqref{eq:ddr.scheme.momentum} with $\uvec{v}_h=\qIcurl{k}{h}(\GRAD\varphi)$ and \eqref{eq:ns.curl.momentum} with $\bvec{v}=\GRAD\varphi \in \bvec{L}^4(\Omega)$ (since $\varphi\in H^2(\Omega)$). The commutation property $\qIcurl{k}{h}(\GRAD\varphi)=\uG{k}{h}(\qIgrad{k}{h}\varphi)$ and the complex properties $\uC{k}{h}\circ\uG{k}{h}=\uvec{0}_h$ and $\CURL\circ\GRAD=0$ then give
\begin{align}
  \term_2={}&\int_\Omega [(\Pdiv{k}{h}\uC{k}{h}\uvec{u}_h)\times\Pcurl{k}{h}\uvec{u}_h]\cdot\Pcurl{k}{h}\qIcurl{k}{h}(\GRAD\varphi)
  -\int_\Omega [(\CURL\bvec{u}\times\bvec{u})]\cdot\GRAD\varphi\nonumber\\
  &-(\Icurl{k}{h}\bvec{f},\qIcurl{k}{h}(\GRAD\varphi))_{\CURL,h}
    +\int_\Omega \bvec{f}\cdot\GRAD\varphi\nonumber\\
  ={}&\term_{2,1}+\term_{2,2}.
  \label{eq:def.T2.T21.T22}
\end{align}
Rearranging the trilinear product as
\[
(\bvec{a}\times\bvec{b})\cdot\bvec{c}-(\bvec{a}'\times\bvec{b}')\cdot\bvec{c}'=
([\bvec{a}-\bvec{a'}]\times\bvec{b})\cdot\bvec{c}+(\bvec{a}'\times[\bvec{b}-\bvec{b'}])\cdot\bvec{c}
+(\bvec{a}'\times\bvec{b}')\cdot[\bvec{c}-\bvec{c}']
\]
and using H\"older inequalities with exponents $(2,4,4)$, we obtain a bound on $\term_{2,1}$:
\begin{align*}
  |\term_{2,1}|&\lesssim \norm{\bvec{L^2}(\Omega)}{\Pdiv{k}{h}\uC{k}{h}\uvec{u}_h-\CURL\bvec{u}}
  \norm{\bvec{L}^4(\Omega)}{\Pcurl{k}{h}\uvec{u}_h}\norm{\bvec{L}^4(\Omega)}{\Pcurl{k}{h}\qIcurl{k}{h}(\GRAD\varphi)}\nonumber\\
  &\quad+\norm{\bvec{L^2}(\Omega)}{\CURL\bvec{u}}
  \norm{\bvec{L}^4(\Omega)}{\Pcurl{k}{h}\uvec{u}_h-\bvec{u}}\norm{\bvec{L}^4(\Omega)}{\Pcurl{k}{h}\qIcurl{k}{h}(\GRAD\varphi)}\nonumber\\
  &\quad+\norm{\bvec{L^2}(\Omega)}{\CURL\bvec{u}}
  \norm{\bvec{L}^4(\Omega)}{\bvec{u}}\norm{\bvec{L}^4(\Omega)}{\Pcurl{k}{h}\qIcurl{k}{h}(\GRAD\varphi)-\GRAD\varphi}.
\end{align*}
By Sobolev embedding and \eqref{eq:reg.varphi}, we have $\norm{\bvec{L}^4(\Omega)}{\GRAD\varphi}\lesssim \norm{\bvec{H}^1(\Omega)}{\GRAD\varphi}\lesssim \norm{L^2(\Omega)}{g}$. Combining this with \eqref{eq:conv.PqI.L4} and the $L^4$-bound on $\Pcurl{k}{h}\uvec{u}_h$ (resulting, e.g., from \eqref{eq:conv.Pcurlh.u}, we infer
\[
|\term_{2,1}|\lesssim \left(\norm{\bvec{L^2}(\Omega)}{\Pdiv{k}{h}\uC{k}{h}\uvec{u}_h-\CURL\bvec{u}}
+\norm{\bvec{L}^4(\Omega)}{\Pcurl{k}{h}\uvec{u}_h-\bvec{u}}
+h^{\frac14}\right)\norm{L^2(\Omega)}{g}.
\]
Turning to $\term_{2,2}$, we write
\begin{align*}
  |\term_{2,2}|&\le\left|\int_\Omega\bvec{f}\cdot\GRAD\varphi-\int_{\Omega}\Pcurl{k}{h}\Icurl{k}{h}\bvec{f}\cdot\GRAD\varphi\right|\\
  &\quad+\left|\int_{\Omega}\Pcurl{k}{h}\Icurl{k}{h}\bvec{f}\cdot\GRAD\varphi-(\Icurl{k}{h}\bvec{f},\qIcurl{k}{h}\GRAD\varphi)_{\CURL,h}\right|\\
  &\lesssim\norm{\bvec{L}^2(\Omega)}{\bvec{f}-\Pcurl{k}{h}\Icurl{k}{h}\bvec{f}}\norm{\bvec{L}^2(\Omega)}{\GRAD\varphi}\\
  &\quad+\norm{\bvec{L}^2(\Omega)}{\GRAD\varphi-\vlproj{k}{h}\GRAD\varphi}\norm{\CURL,h}{\Icurl{k}{h}\bvec{f}},
\end{align*}
where we have used, in the conclusion, the primal consistency of the inner product for $\qIcurl{k}{h}$ stated in \cite[Theorem 8]{Chaumont-Frelet.Droniou.ea:25}. Invoking the approximation properties of $\vlproj{k}{h}$, the estimate \eqref{eq:reg.varphi} on $\varphi$, and the bound \eqref{eq:bound.If.compactness} on $\Icurl{k}{h}\bvec{f}$, we infer
\[
|\term_{2,2}|\lesssim \left(\norm{\bvec{L}^2(\Omega)}{\bvec{f}-\Pcurl{k}{h}\Icurl{k}{h}\bvec{f}}+h\right)\norm{L^2(\Omega)}{g}.
\]
Combining the bounds on $\term_{2,1}$ and $\term_{2,2}$ in \eqref{eq:def.T2.T21.T22} and plugging the result, together with \eqref{eq:conv.pressure.2}, into \eqref{eq:conv.pressure.1}, we see that, since $h\lesssim h^{\frac14}$,
\[
\begin{aligned}
  \int_\Omega (\Pgrad{k+1}{h}\underline{p}_h-p)g
  \lesssim \bigg(
  &h^{\frac14}
  + \norm{\bvec{L^2}(\Omega)}{\Pdiv{k}{h}\uC{k}{h}\uvec{u}_h-\CURL\bvec{u}}
  \\
  &+\norm{\bvec{L}^4(\Omega)}{\Pcurl{k}{h}\uvec{u}_h-\bvec{u}}
  + \norm{\bvec{L}^2(\Omega)}{\bvec{f}-\Pcurl{k}{h}\Icurl{k}{h}\bvec{f}}
  \bigg) \norm{L^2(\Omega)}{g}.
\end{aligned}  
\]

  Polynomial consistency and the pressure normalisation give $\int_\Omega(\Pgrad{k+1}{h}\underline{p}_h - p) = 0$; hence, the $L^2$-norm of $\Pgrad{k+1}{h}\underline{p}_h-p$ can be recovered by duality against zero-average functions.
Taking the supremum over $g\in L^2(\Omega)$ such that $\int_\Omega g=0$ and
$\norm{L^2(\Omega)}{g}=1$ therefore leads to
\begin{multline*}
  \norm{L^2(\Omega)}{\Pgrad{k+1}{h}\underline{p}_h-p}\\
  \lesssim h^{\frac14}+
  \norm{\bvec{L^2}(\Omega)}{\Pdiv{k}{h}\uC{k}{h}\uvec{u}_h-\CURL\bvec{u}}
  +\norm{\bvec{L}^4(\Omega)}{\Pcurl{k}{h}\uvec{u}_h-\bvec{u}}
  +\norm{\bvec{L}^2(\Omega)}{\bvec{f}-\Pcurl{k}{h}\Icurl{k}{h}\bvec{f}}
\end{multline*}
and the convergence properties \eqref{eq:conv.Ch.u}, \eqref{eq:conv.Pcurlh.u} and \eqref{eq:conv.PI.Linfty} complete the proof of \eqref{eq:conv.Pgradh.p}.


\section{Numerical tests}
\label{sec:num.tests}

We consider the classical lid-driven cavity test on the unit cube $\Omega \coloneqq (0,1)^3$.
  Two sets of boundary conditions are considered:
  the first one consists of enforcing the tangential velocity on the lid and homogeneous natural boundary conditions on the other faces;
  the second one involves additional penalty terms in the spirit of \cite{Wietse-M.Wouter.ea:26} to incorporate no-slip conditions, thus mimicking more closely the usual formulation of the problem when the viscous terms are written in grad--grad form.
  As we will see, the two sets lead to entirely different solutions.

\subsection{Mixed boundary conditions}
\label{sec:num.tests:mixed}

We first consider
\begin{itemize}
\item Essential boundary conditions $p = 0$ and $\bvec{n} \times (\bvec{u} \times \bvec{n}) =  (0,1,0)$ on the lid;
\item Homogeneous natural boundary conditions $(\CURL \bvec{u} )\times \bvec{n} = 0$ and $\bvec{u} \cdot \bvec{n} = 0$ on the remaining boundary faces.
\end{itemize}
The essential boundary condition on the lid does not prevent inflow or outflow,
and the natural boundary condition does not enforce a no-slip condition.
This allows, in particular, tangential flow along the walls.

\begin{figure}
  \centering
  \begin{subfigure}[b]{0.48\textwidth}
    \centering
    \includegraphics[width=\textwidth]{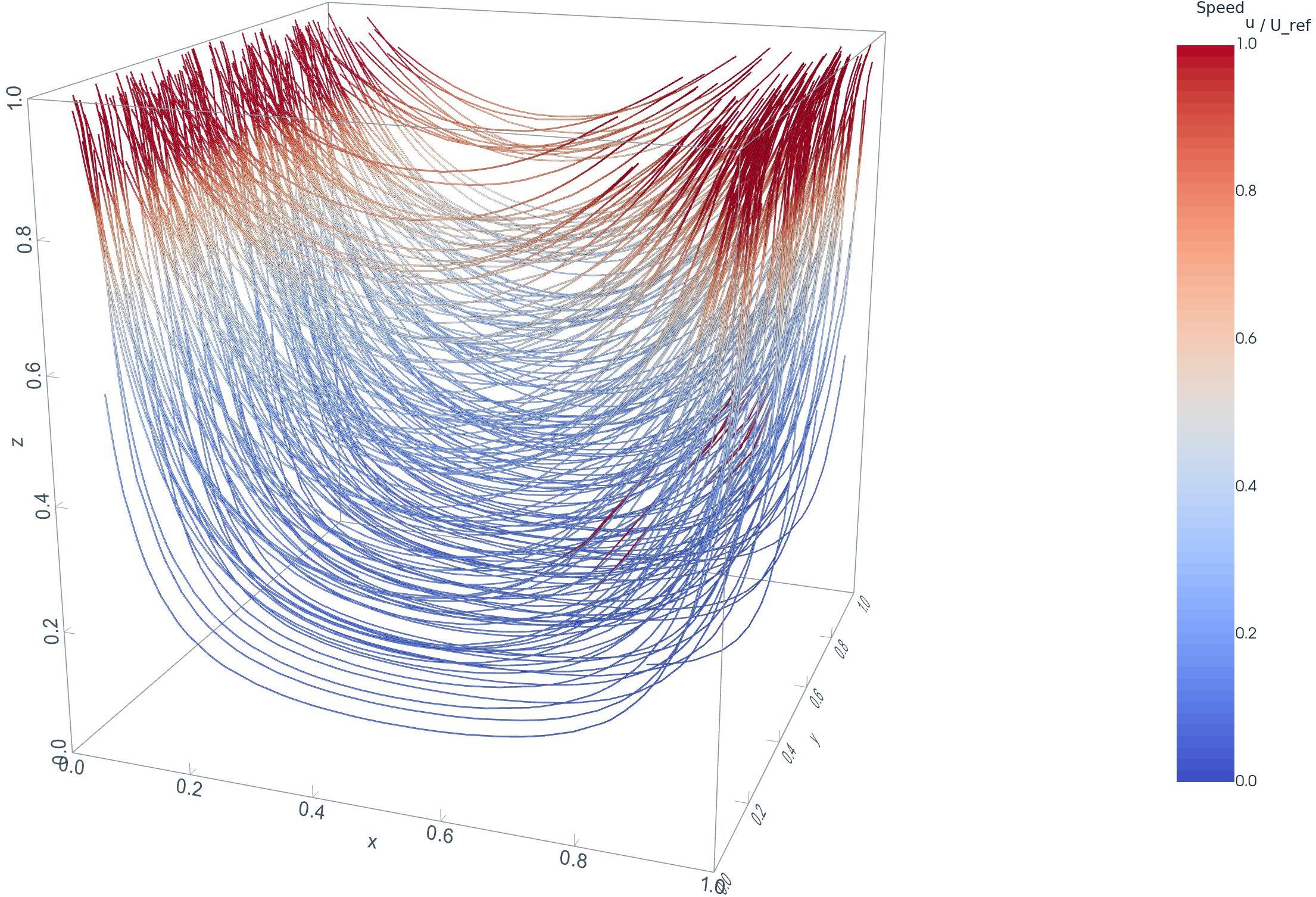}
    \caption{Reynolds number $400$.}
  \end{subfigure}
  \hfill
  \begin{subfigure}[b]{0.48\textwidth}
    \centering
    \includegraphics[width=\textwidth]{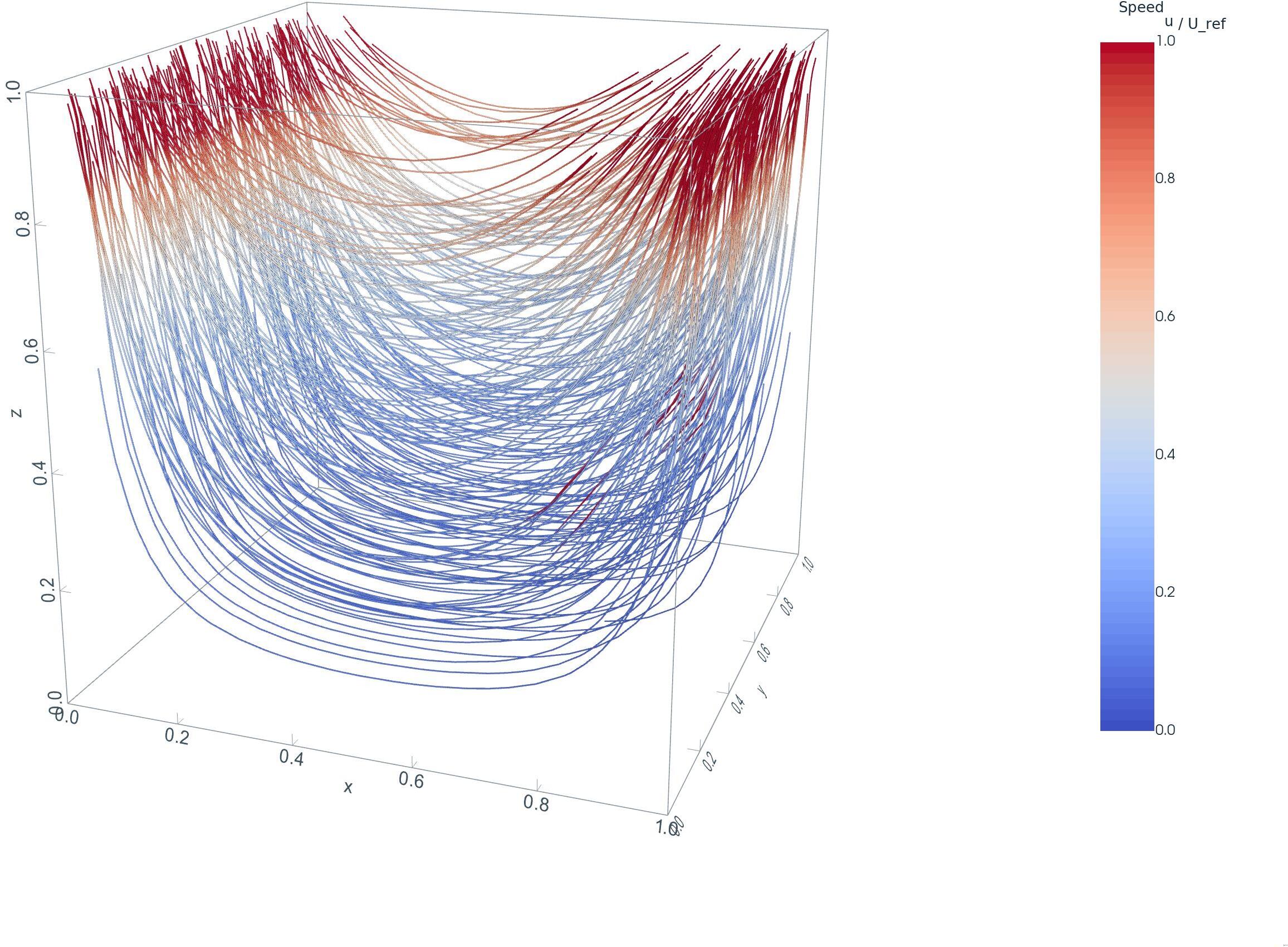}
    \caption{Reynolds number $10000$.}
  \end{subfigure}
  \caption{Velocity streamlines for the lid-driven cavity.}
  \label{fig:streamlineCavity}
\end{figure}

  The numerical scheme was implemented in the \hardcore platform\footnote{\url{https://github.com/jdroniou/HArDCore}}.
The meshes for these tests are Voronoi meshes, and we consider the serendipity DDR complex with
an element-specific choice of serendipity spaces using the procedure described in \cite{Di-Pietro.Droniou:23*1}, and polynomial
degrees $k = 0,1,2$.
The nonlinear problem was solved using Newton's method, with pseudo-time iterations and gradually increased Reynolds number to its target value, to ensure convergence at large Reynolds numbers.
  The linear systems were solved using the SparseLU direct solver, with initial relaxation parameter of $0.5$, which was gradually increased after successful Newton iterations.
  The initial system with Reynolds number $1$ was solved in a few Newton iterations.
  The subsequent systems, obtained by increasing the Reynolds number, only required a single Newton iteration each.
We computed the solution for Reynolds numbers $100$, $400$, $1000$, and $10000$.
The velocity streamlines for Reynolds numbers $400$ and $10000$ computed using a polynomial degree $k=1$ are shown in Figure \ref{fig:streamlineCavity}.
In contrast with the solution obtained using no-slip boundary conditions, we observed little change in the solution as the Reynolds number is increased.

\FloatBarrier

\subsection{Weakly enforced no-slip boundary conditions}
\label{sec:num.tests:no-slip}

To obtain a problem that more closely resembles the standard three-dimensional lid-driven cavity,
we adapted the method introduced in \cite{Wietse-M.Wouter.ea:26} for the Stokes problem.
As demonstrated therein, enforcing no-slip conditions as essential boundary conditions leads to ill-posed problems on certain meshes.
This difficulty can be circumvented by using Nitsche-type techniques to enforce this condition weakly.
In the present context, this amounts to adding the following boundary terms to the discrete momentum conservation equation:
\[
\nu\int_{\partial\Omega} (\bvec{n} \times \Pdiv{k}{h} \uC{k}{h} \uvec{u}_h ) \cdot \trFt \uvec{v}_h
- \nu \int_{\partial \Omega} \trFt \uvec{u}_h \cdot (\bvec{n} \times \Pdiv{k}{h} \uC{k}{h} \uvec{v}_h )
+ \nu \int_{\partial \Omega} \frac{1}{h} \trFt \uvec{u}_h \cdot \trFt \uvec{v}_h.
\]
The resulting boundary conditions are:
\begin{itemize}
  \item Lid-driven velocity $\bvec{u} = (1,0,0)$ on the lid;
  \item No-slip condition $\bvec{u} = 0$ on the remaining boundary faces.
\end{itemize}

  The corresponding numerical scheme, along with a standard HHO formulation \cite{Botti.Di-Pietro.ea:19}, was implemented in the \polyrust library\footnote{\url{https://polyrust.org}} \cite{Di-Pietro:26}.
  The versions of both schemes corresponding to $k = 1$ were considered on a non-uniform Cartesian orthogonal mesh obtained using Gmsh\footnote{\url{https://gmsh.info}}, with tapered element sizes ranging from $0.1$ at the centre to $1.25 \times 10^{-2}$ towards the boundary.
  In both cases, convergence at higher Reynolds numbers was obtained using a continuation procedure, increasing the Reynolds number from $\Reynolds = 100$ to $\Reynolds = 250$ in increments of $25$.
  On this mesh, convective effects became more pronounced around $\Reynolds = 300$, preventing the iterative solver from converging to the desired tolerance.
  Simulations at higher Reynolds numbers would therefore have required finer meshes and/or the use of convective stabilisation.
  Stabilising convection for curl--curl formulations requires further investigation; this topic is left for future work.

  For each value of the Reynolds number, the nonlinear algebraic problem was solved with a quasi-Newton method incorporating an adaptive stopping criterion for the iterative linear solver inspired by \cite{Eisenstat.Walker:96} with the following parameters: initial value $0.1$, maximum value $0.5$, scaling parameter $0.9$, and exponent $1.5$.
  The nonlinear residual tolerance was set to $10^{-8}$.

  At each Newton iteration, the statically condensed linear systems were solved using PETSc with right preconditioning and modified Gram--Schmidt orthogonalisation.
  The following additional settings were identified by comparing a number of potentially suitable configurations.
  For DDR, flexible GMRES with restart length 500 was combined with a restricted additive Schwarz preconditioner with zero overlap and local sparse LU factorisations.
  METIS nested-dissection ordering was applied both globally and within the subdomains, with a diagonal-shift parameter of $10^{-12}$ to avoid zero pivots.
  For HHO, GMRES with restart length 100 was combined with restricted additive Schwarz with one layer of overlap and local incomplete LU factorisations with zero fill-in, ILU(0), using natural ordering.
  In both cases, the preconditioners were rebuilt for each Newton Jacobian.

\begin{figure}
  \centering
  \begin{subfigure}[b]{0.48\textwidth}
    \centering
    \includegraphics[width=\textwidth]{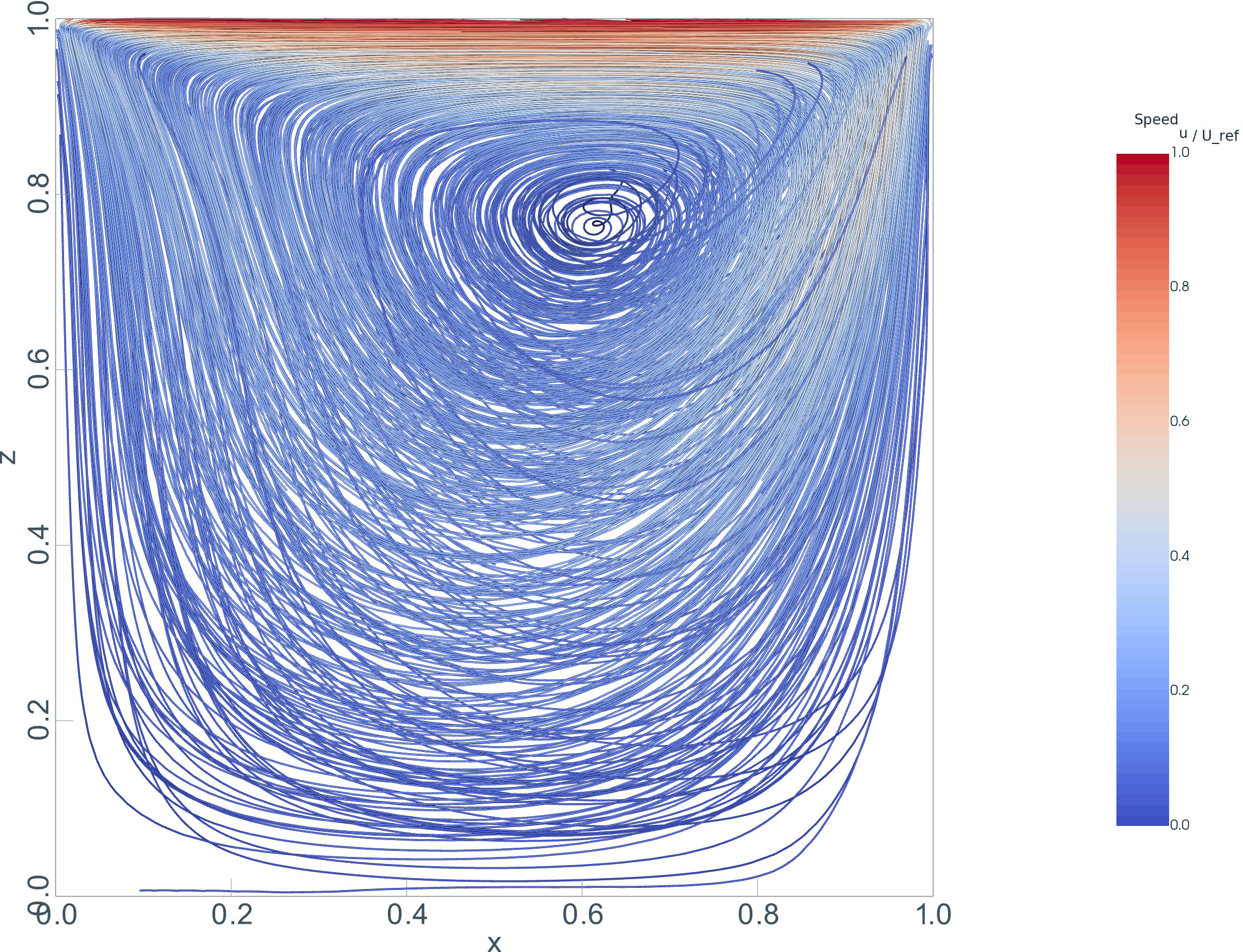}
  \end{subfigure}
  \hfill
  \begin{subfigure}[b]{0.48\textwidth}
    \centering
    \includegraphics[width=\textwidth]{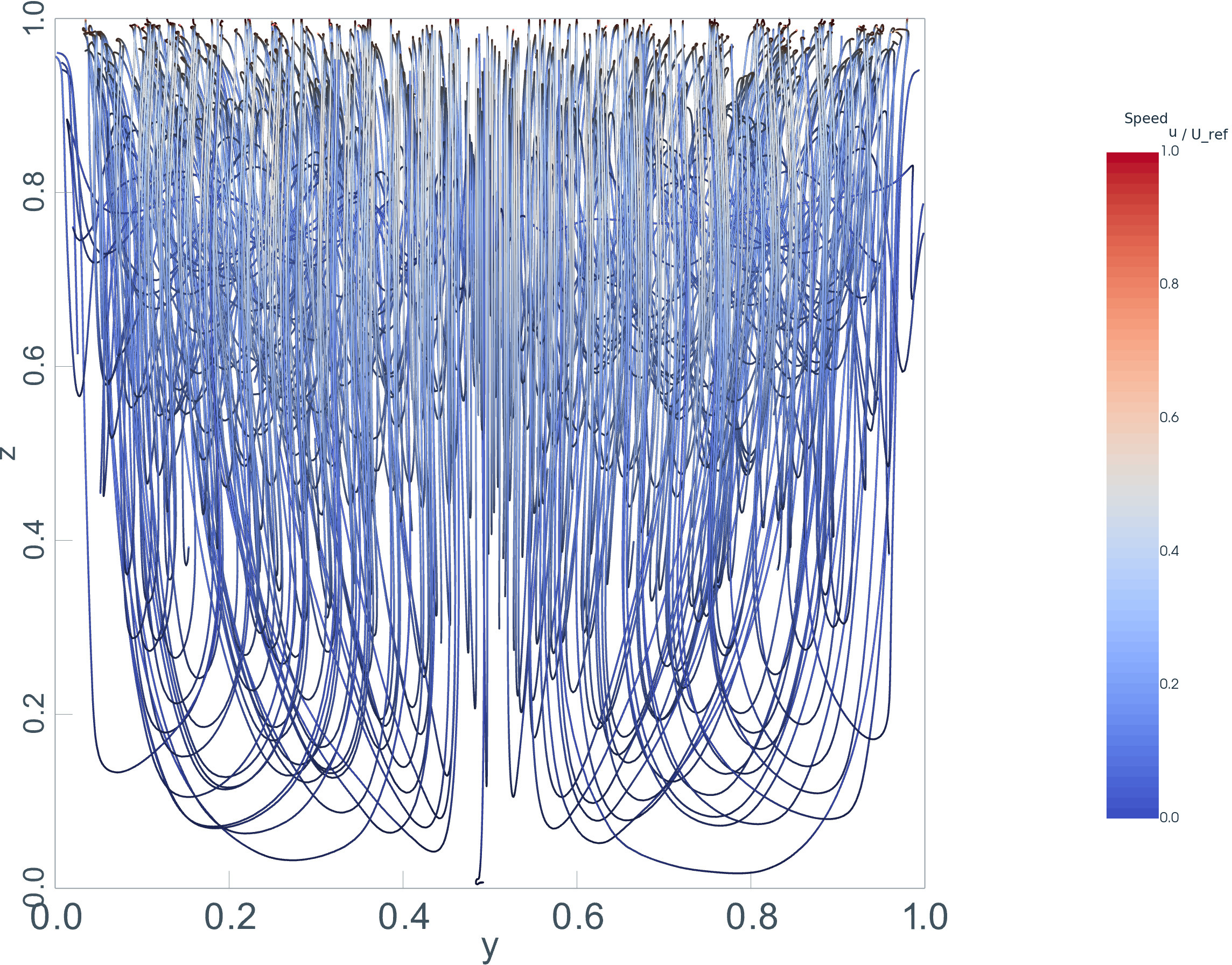}
  \end{subfigure}
  \caption{Side and front views of the velocity streamlines computed using the DDR method for the lid-driven cavity with the weakly enforced no-slip conditions of Section \ref{sec:num.tests:no-slip} at Reynolds number $100$.}
  \label{fig:streamlines-ddr}
\end{figure}

\begin{figure}
  \begin{subfigure}[t]{0.95\textwidth}
    \centering
    \includegraphics[width=7.0cm]{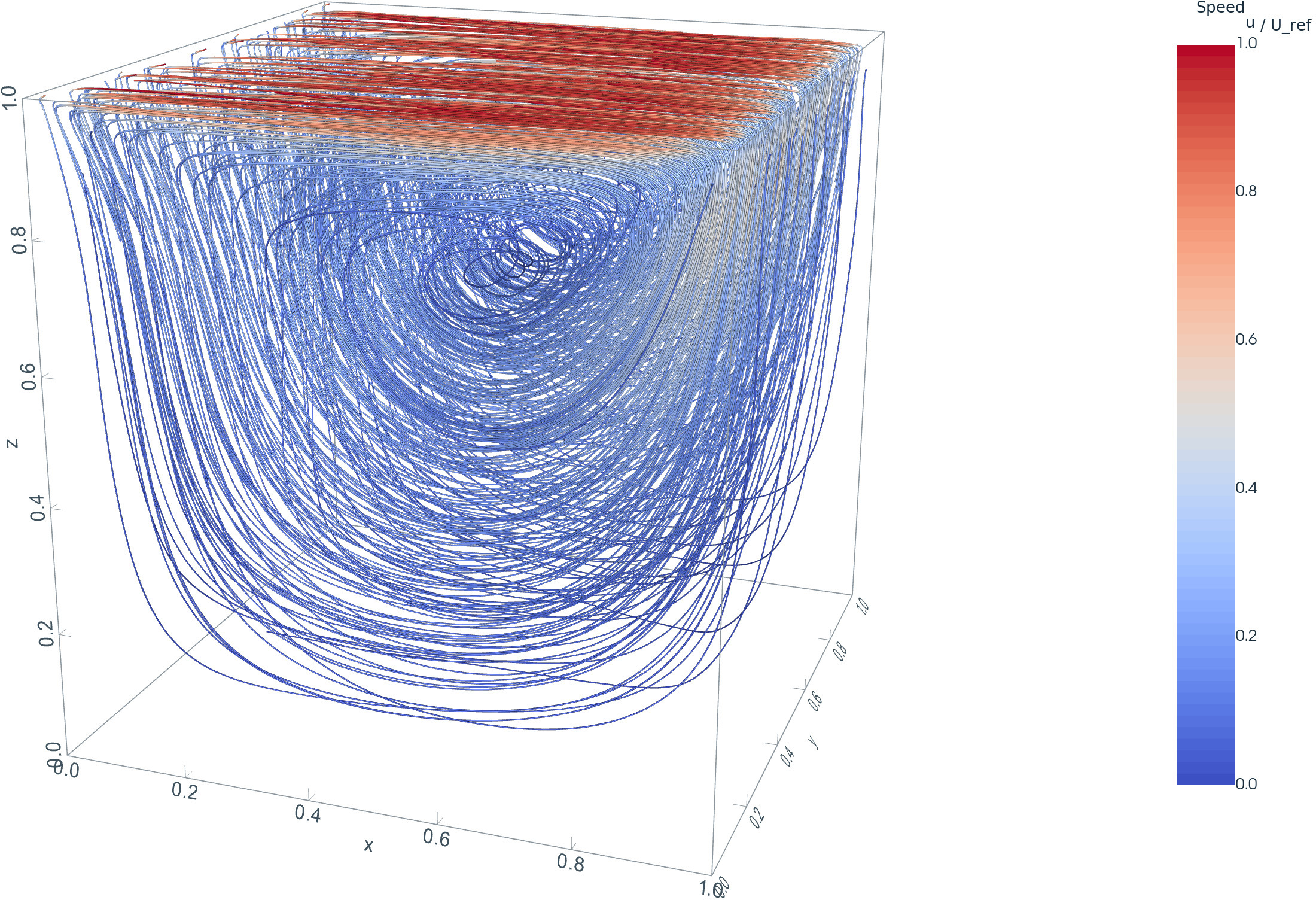}
    \hspace{0.5cm}
    \includegraphics[width=7.0cm]{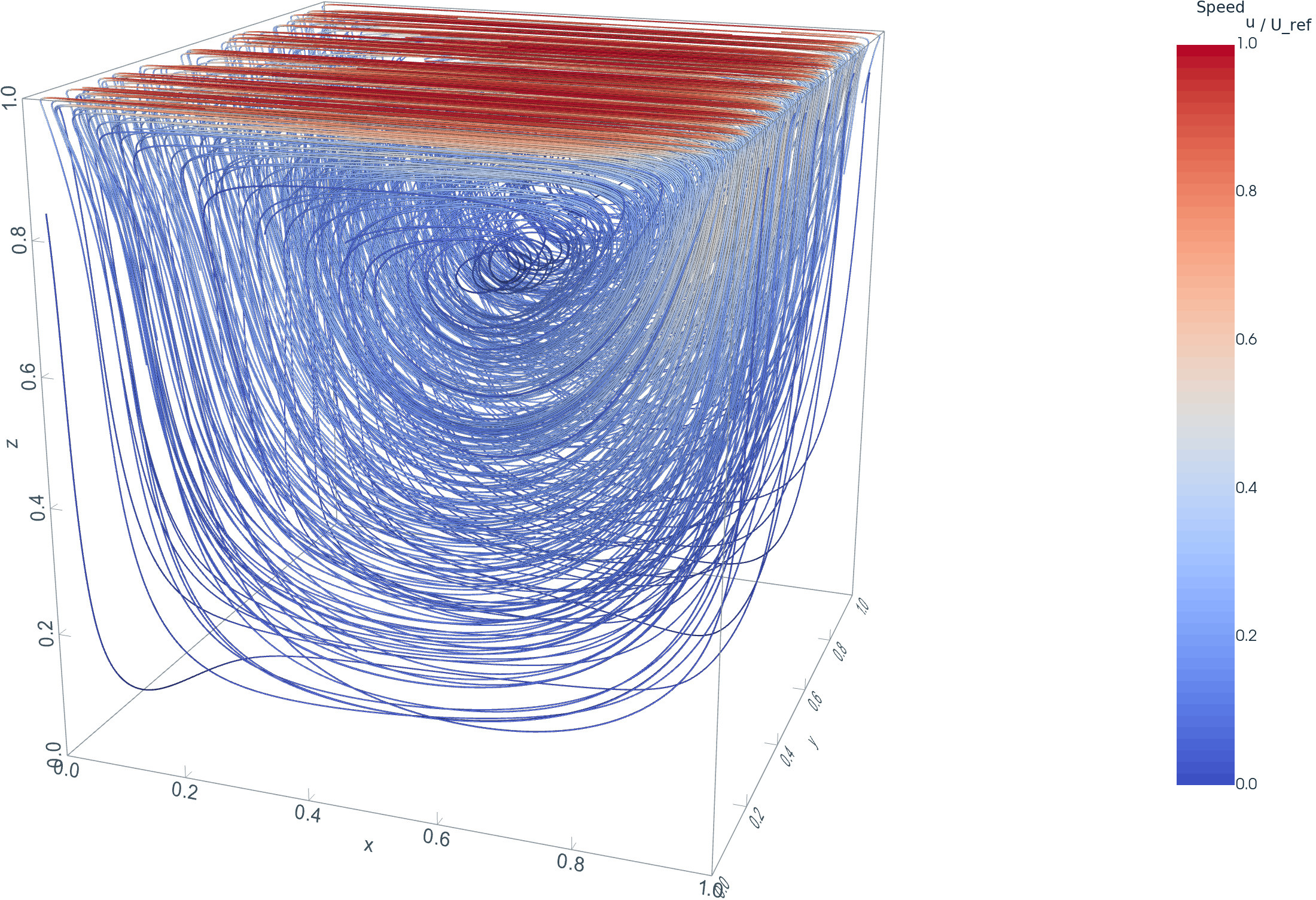}
    \caption{$\Reynolds = 100$}
  \end{subfigure}
  \vspace{0.25cm} \\
  \begin{subfigure}[t]{0.95\textwidth}
    \centering
    \includegraphics[width=7.0cm]{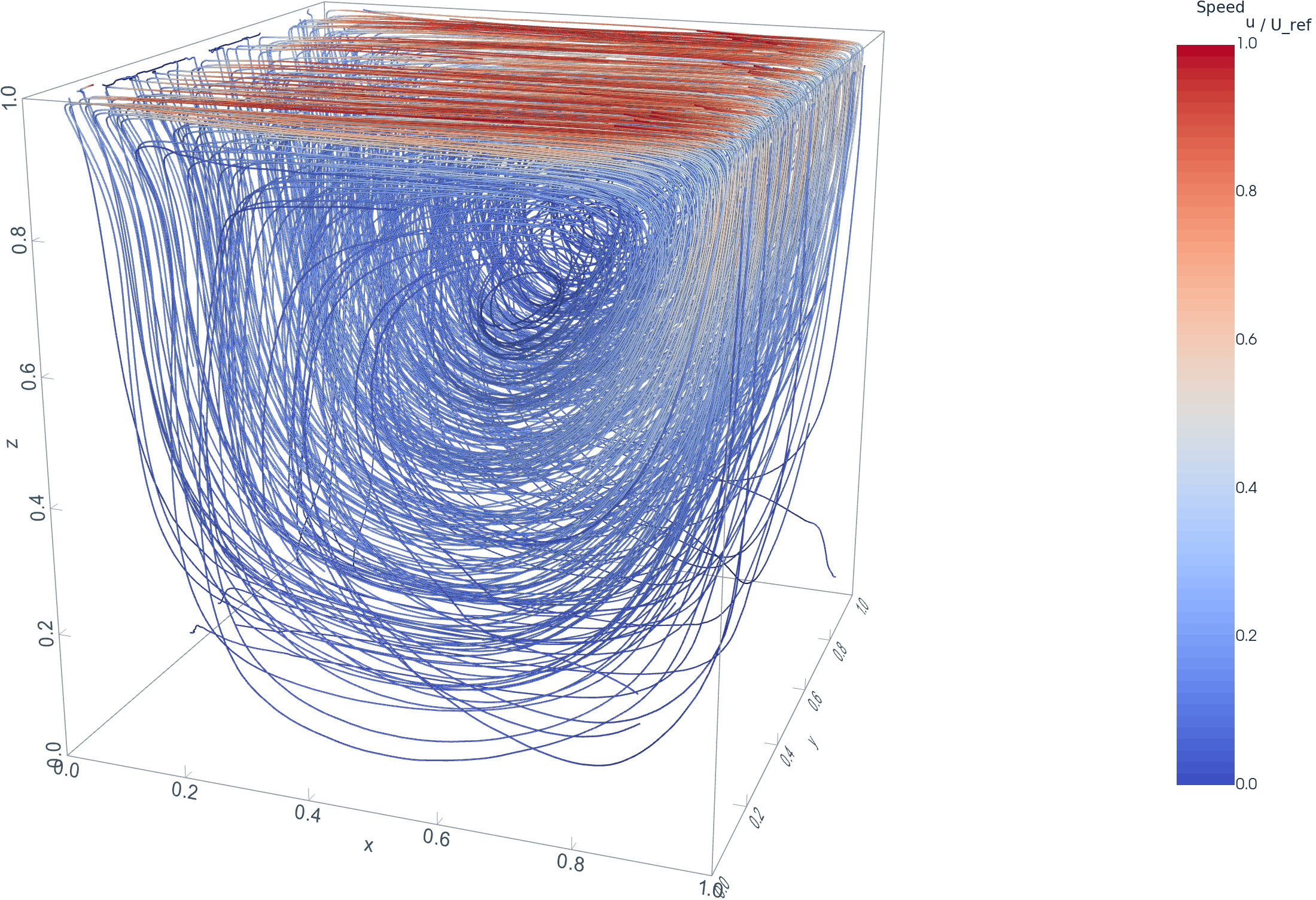}
    \hspace{0.5cm}
    \includegraphics[width=7.0cm]{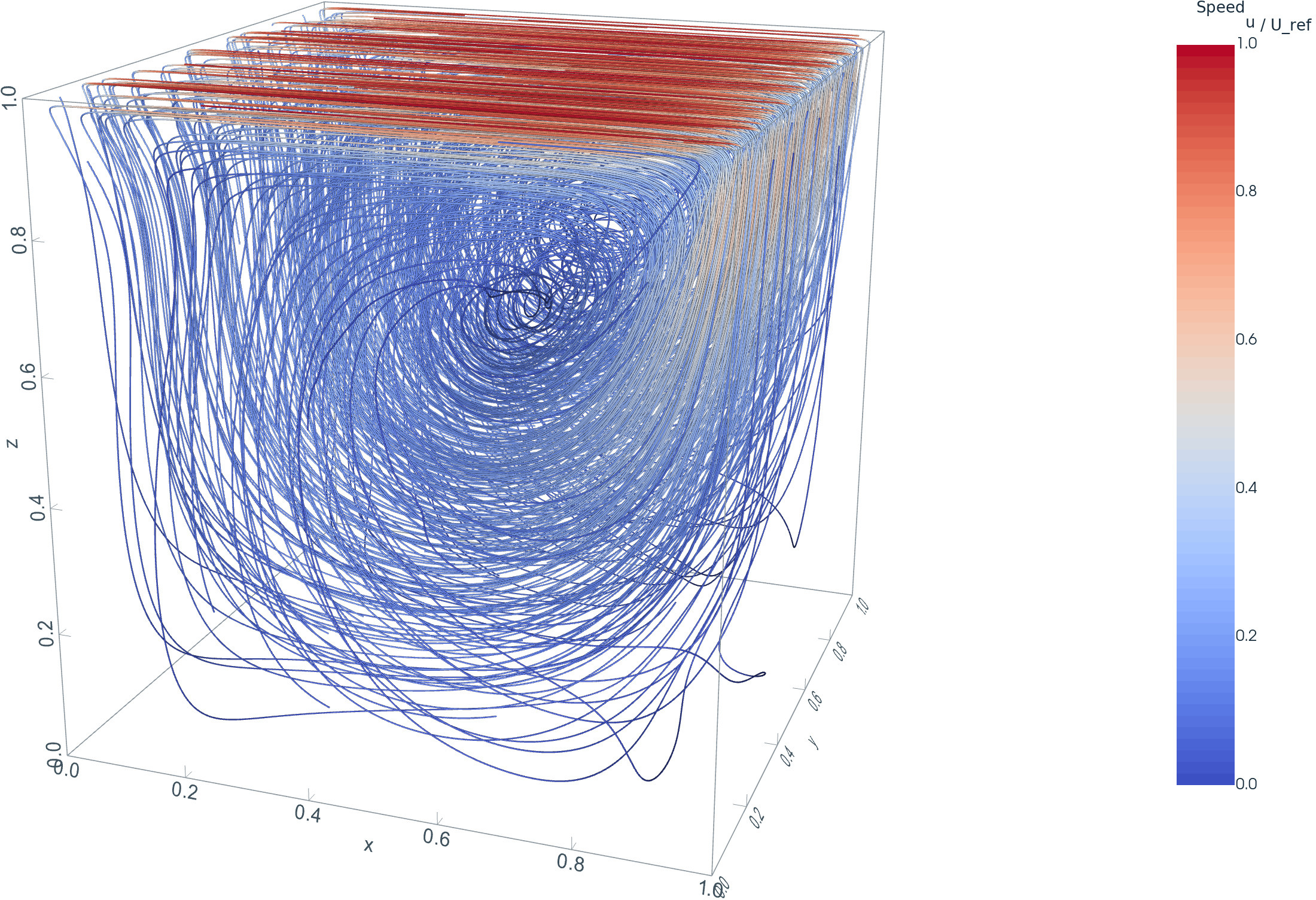}
    \caption{$\Reynolds = 250$}
  \end{subfigure}
  \caption{Comparison between the present DDR scheme (left) and the HHO solution (right) for the three-dimensional lid-driven cavity problem.
    For both schemes, we used $k = 1$ and a hexahedral mesh refined towards the corners.
  }
  \label{fig:cavity:ddr-hho-comparison}
\end{figure}

\begin{figure}
  \centering
  \includegraphics[width=0.4\textwidth]{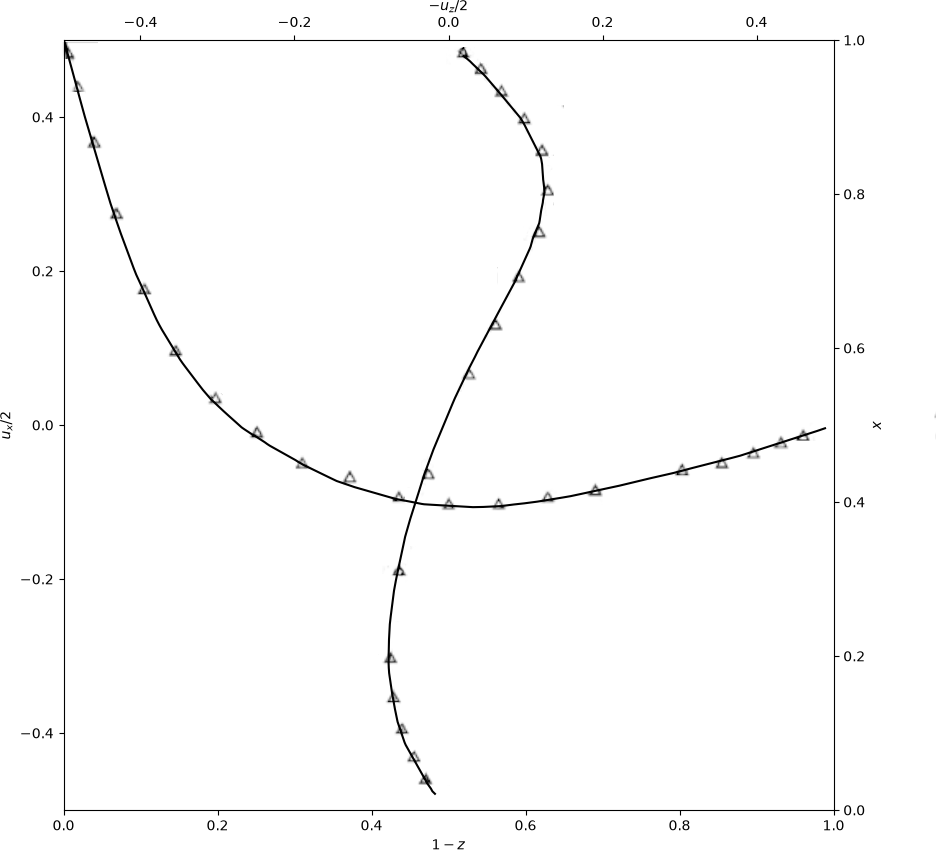}
  \caption{Fluid velocity along the centrelines $(0.5,0.5,z)$ and $(x,0.5,0.5)$.
    The triangles represent the results of \cite{Hwar-C.Richard-S.ea:87}, extracted from their figures.}
  \label{fig:cmp_profiles}
\end{figure}

Figure \ref{fig:streamlines-ddr} shows side and front views of the velocity streamlines at Reynolds number $\Reynolds = 100$.
The Nitsche-type method correctly enforces both the normal and tangential components of the velocity fields.
Compared with the mixed boundary conditions considered in Section \ref{sec:num.tests:mixed}, this configuration produces vortices.
This is illustrated in Figure \ref{fig:cavity:ddr-hho-comparison}, which shows side by side the DDR and HHO solutions on the same meshes for $\Reynolds \in \{ 100, 250\}$.
The streamlines were drawn using the same seeds and parameters, so a direct comparison is possible.
The flow exhibits a distinctly three-dimensional feature: an inward spiralling motion perpendicular to the main vortex.
Both solutions exhibit essentially the same behaviour, with only minor differences towards the boundary.
The results agree well with previous studies, as illustrated by the comparison with \cite{Hwar-C.Richard-S.ea:87} in Figure \ref{fig:cmp_profiles}.


\section*{Acknowledgements}

Funded by the European Union (ERC Synergy, NEMESIS, project number 101115663).
Views and opinions expressed are however those of the authors only and do not necessarily reflect those of the European Union or the European Research Council Executive Agency. Neither the European Union nor the granting authority can be held responsible for them.


\printbibliography

\end{document}